\documentclass[a4paper,12pt]{amsart}

\usepackage[french]{babel}
\usepackage{amsmath, enumerate, amsfonts, amssymb, amsthm}

\newtheorem{defin}{\bf D\'ef\mbox{}inition}[section]
\newtheorem{theo}[defin]{\bf Th\'eor\`eme}
\newtheorem{prop}[defin]{\bf Proposition}
\newtheorem{lem}[defin]{\bf Lemme}
\newtheorem{cor}[defin]{\bf Corollaire}

\newtheorem{rem}[defin]{\bf Remarque}

\newtheorem*{aff*}{\bf Affirmation}
\newtheorem*{rem*}{\bf Remarque}
\newtheorem*{nota*}{\bf Notation}
\newtheorem*{prop*}{\bf Proposition}

\newtheorem{num}[defin]{}

\newtheorem{theoI}{\bf Th\'eor\`eme}
\newtheorem{propI}{\bf Proposition}
\newtheorem{corI}{\bf Corollaire}
\newtheorem{remI}{\bf Remarque}

\newcommand{\dps}{\displaystyle}

\newcommand{\C}{\mathbb{C}}
\newcommand{\N}{\mathbb{N}}
\newcommand{\Z}{\mathbb{Z}}
\newcommand{\Q}{\mathbb{Q}}
\renewcommand{\k}{\mathbf{k}}

\newcommand{\ring}{\mathcal{R}}

\newcommand{\CC}{\mathcal{C}}
\newcommand{\FF}{\mathcal{F}}

\newcommand{\QQ}{\mathcal{Q}}
\newcommand{\modQ}{{\mathrm{mod}\QQ}}
\newcommand{\PP}{\mathcal{P}}
\newcommand{\II}{\mathcal{I}}
\newcommand{\GG}{\mathcal{G}}

\newcommand{\spec}{\mathrm{Spec}}
\newcommand{\specm}{\mathrm{Specm}}
\newcommand{\Frac}{\mathrm{Frac}}

\newcommand{\An} {\mathbf{A}_n}

\newcommand{\sdt}{\langle s, \partial_t \rangle}

\renewcommand{\O}{\mathcal{O}} 

\newcommand{\D}{\mathcal{D}} 

\newcommand{\FDn}{\hat{\mathcal{D}}_n}
\newcommand{\FDnpun}{\hat{\mathcal{D}}_{n+1}}

\newcommand{\dxi}{\partial _{x_i}}
\newcommand{\dx}[1]{\partial _{x_{#1}}}
\newcommand{\ddx}{\partial _x}
\newcommand{\ddt}{\partial _t}

\newcommand{\dxsur}[2]{\frac{\partial {#1}}{\partial x_{#2}}}

\newcommand{\B}{\mathcal{B}}

\newcommand{\ord}{\mathrm{ord}}
\newcommand{\cp}{\mathrm{cp}} 
\newcommand{\tp}{\mathrm{tp}} 
\newcommand{\DN}{\mathcal{N}} 
\newcommand{\Exp}{\mathrm{Exp}} 

\title[$b$-fonction g\'en\'erique locale]{Polyn\^ome de Bernstein-Sato
  g\'en\'erique local}

\author{Rouchdi BAHLOUL}
\address{Laboratoire de Math\'ematiques,
Universit\'e de Versailles Saint-Quentin-en-Yvelines,
45 avenue des \'Etats-Unis, 78035 Versailles, France}
\email{bahloul@math.uvsq.fr}

\begin{document}

\begin{abstract}
\'Etant donn\'ee une famille de fonctions analytiques en $0\in \C^n$
param\'etr\'ee par un espace lisse, nous \'etudions le polyn\^ome de
Bernstein de la fibre sur une vari\'et\'e irr\'eductible $V$ de l'espace
des param\`etres et nous montrons qu'il est g\'en\'eriquement constant.
Nous montrons que ce polyn\^ome $b$ satisfait une \'equation fonctionnelle
g\'en\'erique sur $V$ et l'on d\'erive une stratification constructible
de l'espace des param\`etres par le polyn\^ome de Bernstein de la fibre.
Lorsque l'hypersurface admet g\'en\'eriquement une singularit\'e unique
en $0\in \C^n$ nous montrons que $b$
est le polyn\^ome de Bernstein g\'en\'erique au sens de
Brian\c{c}on-Geandier-Maisonobe.
Les outils utilis\'es sont une
g\'en\'eralisation formelle d'un algorithme de Oaku calculant le polyn\^ome
de Bernstein local et les bases standard g\'en\'eriques r\'ecemment
\'etudi\'ees par l'auteur.
\end{abstract}

\subjclass[2000]{Primary 32S30; Secondary 16S32, 13P99}
\keywords{Polyn\^ome de Bernstein-Sato, d\'eformations de singularit\'es,
bases standard param\'etriques}

\maketitle


\section*{Introduction et motivations}\label{sec:intro}

Le polyn\^ome de Bernstein (ou $b$-fonction) a \'et\'e introduit de mani\`ere
ind\'ependante par I.~N. Bernstein \cite{bernstein} et M. Sato \cite{sato}.
Son existence a \'et\'e d\'emontr\'ee par Bernstein dans le cas polynomial et
par J.~E.~Bj\"ork \cite{bjorkP} dans le cas analytique et formel (voir
aussi \cite{bjorkB}) ainsi que par M. Kashiwara \cite{kashiwara} qui
d\'emontra en plus la rationalit\'e des racines du polyn\^ome de Bernstein
analytique. Dans le cas de plusieurs fonctions analytiques, l'existence
de polyn\^omes de Bernstein revient \`a C.~Sabbah (\cite{sabbah1, sabbah2},
voir aussi \cite{gyoja} et \cite{Comp}). Ici, nous nous int\'eressons au
polyn\^ome de Bernstein d'une fonction analytique d\'ependant de param\`etres.

On sait depuis les travaux de D.~T. L\^e et C.~P. Ramanujam \cite{le, le-r}
qu'une d\'eformation \`a nombre de Milnor $\mu$ constant
d'une singularit\'e isol\'ee d'hypersurface conserve son type topologique
ainsi que la classe de conjugaison de sa monodromie locale.
Ceci combin\'e aux travaux de B. Malgrange \cite{malgrange} liant la
monodromie et les racines du polyn\^ome de Bernstein local nous dit
que les racines de ce dernier restent inchang\'ees modulo $\Z$.
Cependant ses racines ne sont pas constantes (voir par exemple T. Yano
\cite{yano})

F. Geandier \cite{geandierCRAS, geandierComp} a \'etudi\'e de
mani\`ere \'etendue le polyn\^ome de Bernstein associ\'e \`a une
d\'eformation \`a un param\`etre d'une hypersurface \`a singularit\'e
isol\'ee avec une \'etude du polyn\^ome de Bernstein g\'en\'erique et
relatif (ou ``en famille''). J.~Brian\c{c}on, F.~Geandier et Ph.~Maisonobe
\cite{bgm} ont g\'en\'eralis\'e et compl\'et\'e l'\'etude pr\'ec\'edente
au cas de plusieurs param\`etres, toujours dans le cas d'une d\'eformation
d'une singularit\'e isol\'ee.

Dans \cite{bgmm}, J. Brian\c{c}on, M. Granger, Ph. Maisonobe et M. Miniconi
ont donn\'e un algorithme de calcul du polyn\^ome de Bernstein pour une
fonction semi-quasi-homog\`ene ou non d\'eg\'en\'er\'ee au sens de
Kouchnirenko. C'est aussi l\`a que la notion de polyn\^ome de Bernstein
g\'en\'erique fut introduite avec le calcul exact dans le cas d'une
d\'eformation semi-universelle \`a deux variables. Dans le m\^eme
esprit, citons \'egalement les travaux de P. Cassou-Nogu\`es
\cite{cassou86, cassou87, cassou88}.

T.~Oaku \cite{oaku} a donn\'e un algorithme (sans conditions)
de calcul du polyn\^ome de Bernstein local et global associ\'e \`a un
polyn\^ome. Ces algorithmes sont bas\'es sur les bases de Gr\"obner
dans des anneaux d'op\'erateurs diff\'erentiels polynomiaux.
Dans \cite{oakuJPAA}, il a initi\'e une \'etude param\'etrique
de ses algorithmes.
Ceci a permis \`a A.~Leykin \cite{leykin} d'obtenir un r\'esultat de
constructibilit\'e concernant le polyn\^ome de Bernstein global pour un
polyn\^ome d\'ependant de param\`etres (voir aussi \cite{bmai02} et
\cite{these}).

Ce rappel historique (non exhaustif) \'etant fait, introduisons le
pr\'esent travail.

On fixe deux entiers $n, m>0$. Pour commencer, consid\'erons une fonction
polynomiale $f=f(x,y) \in \k[x,y]$ o\`u $x=(x_1,\ldots,x_n)$ est le syst\`eme
de variables principales et $y=(y_1,\ldots,y_m)$ est vu comme param\`etre
(et $\k$ est un corps de caract\'eristique z\'ero). Il est bien connu que le
polyn\^ome de Bernstein g\'en\'erique global est non nul (voir Biosca
\cite{biosca, bioscaT} pour $\k=\C$, \cite{PJA} en g\'en\'eral).
De plus nous savons qu'il est \'egal au polyn\^ome de Bernstein usuel
de $f$ vu dans $\mathrm{Frac}(\k[y])[x]$ et nous savons enfin
qu'il est \'egal au polyn\^ome de Bernstein de la fibre g\'en\'erique (voir
\cite{leykin}, voir aussi \cite{bmai02} et \cite{these}).

Plus g\'en\'eralement on peut consid\'erer le polyn\^ome de Bernstein
de la classe $(f)_\QQ$ de $f$ modulo un id\'eal premier $\QQ \subset \k[y]$,
vue dans $\mathrm{Frac}(\k[y]/\QQ)[x]$ (ce que, dans \cite{PJA}, nous avons
appel\'e polyn\^ome de Bernstein g\'en\'erique de $f$ sur $V(\QQ) \subset
\spec(\k[y])$) et l'on montre alors que d'une part il satisfait une
\'equation fonctionnelle ``g\'en\'erique'' (\cite{PJA}) et
d'autre part c'est le polyn\^ome de Bernstein de la fibre g\'en\'erique
sur $V(\QQ)$ (voir les trois r\'ef\'erences ci-dessus).

Maintenant si $f$ est un germe de fonction analytique dans $\C\{x,y\}$,
nous pouvons faire une construction similaire: voir $f$ dans $\CC[[x]]$
avec $\CC=\C\{y\}$ et consid\'erer le polyn\^ome de Bernstein formel de $f$
vu dans $\mathrm{Frac}(\CC)[[x]]$ (dont l'existence est assur\'ee par
\cite{bjorkP}). Plus g\'en\'eralement, \'etant donn\'e $\QQ \in \spec(\CC)$,
on voit $f$ dans $\CC[[x]]$, on consid\`ere sa classe modulo $\QQ$ que
l'on voit dans $\mathrm{Frac}(\CC/\QQ)[[x]]$ et enfin on prend son polyn\^ome
de Bernstein formel $b$.
Avec ce polyn\^ome $b$, a-t-on des r\'esultats similaires \`a ceux du cas
global rappel\'es ci-dessus~? Nous savons que les choses sont plus complexes
dans le cas local puisque par exemple, le polyn\^ome de Bernstein
g\'en\'erique n'existe pas toujours (Biosca \cite{biosca, bioscaT}).
Le but du pr\'esent travail est de mieux comprendre le r\^ole jou\'e par $b$.
Bien qu'il soit d\'efini de mani\`ere alg\'ebrique, nous montrons qu'il a
un r\^ole g\'eom\'etrique naturel et nous faisons le lien avec le travail
de Brian\c{c}on et al. \cite{bgm}.

\noindent
{\bf Note.}
Dans la suite, pour plus de g\'en\'eralit\'e et aussi pour adh\'erer
aux notations g\'en\'eralement utilis\'ees, nous travaillerons sur
un polydisque compact $Z=X\times Y \subset \C^n \times \C^m$ ($X$
pouvant \^etre nul). Si $\O$ d\'esigne le faisceau des
fonctions holomorphes sur $\C^{n+m}$ alors on notera $\O_Z$ les sections
globales du faisceau $\O_{|Z}$ restreint \`a $Z$. L'anneau $\O_Z$ est
noeth\'erien (Frisch, \cite{frisch}). Cette noeth\'erianit\'e nous sera
n\'ecessaire pour terminer la preuve du th\'eor\`eme principal.
Dans la suite, nous n'aurons pas besoin de r\'eduire le diam\`etre des
polydisques sauf celui de $Y$ dans la d\'emonstration du corollaire
\ref{c2} et de la remarque \ref{rem:theo}(c) (nous n'en ferons pas
mention explicite).

\noindent
{\bf Remerciements.} Je remercie Michel Granger qui, en fevrier 2004,
m'a aid\'e \`a trouver une erreur dans une version pr\'eliminaire
ainsi que pour des discussions \'eclairantes.
Ce travail est effectu\'e dans le cadre d'une bourse post-doctorale
FY2003 de la JSPS.

\section{\'Enonc\'e des r\'esultats principaux}

Soient $X \subset \C^n$ et $Y \subset \C^m$ des polydisques compactes
centr\'es en $0$, $Z=X \times Y$ et $f$ une fonction analytique sur $Z$
telle que l'hypersurface $W=f^{-1}(0) \subset Z$ contienne $0$. Pour
$y\in Y$, $W_y$ est l'hypersurface de $X$ d\'efinie par $f_y:x \mapsto
f(x,y)$.

$\D_{Z/Y}$ d\'esigne l'anneau des op\'erateurs diff\'erentiels relatifs.
C'est le sous anneau de $\D_Z$ constitu\'e d'op\'erateurs sans d\'erivation
par rapport aux $y_i$. Suivant \cite{bmai02}, introduisons
$\C\sdt$ l'alg\`ebre $\C[s,\ddt]$ modulo
la relation $\ddt s= s \ddt -s$ et $\D_Z\sdt = \D_Z \otimes \C\sdt$.
Si $t$ d\'esigne une nouvelle variable alors l'identification $s=-\ddt t$
fournit les inclusions d'anneaux~: $\D_Z[s] \subset \D_Z\sdt \subset
\D_{Z \times \C}$.
Cette identification provient du fait que le module libre
$\O_Z[1/f, s] \cdot f^s$ est un $\D_{Z \times \C}$-module et que
l'action de $s$ co\"{\i}ncide avec celle de $- \ddt t$ (voir Malgrange
\cite{malgrange}).

Soit $\CC=\O_Y$ (l'anneau des param\`etres) et $\QQ \in \spec(\CC)$.
On voit $f$ dans $\CC[[x]]$ et on note $[f]_\QQ \in (\CC/\QQ)[[x]]$ la
s\'erie obtenue en prenant la classe modulo $\QQ$ des coefficients de
$f$. Enfin, on note $(f)_\QQ$ la s\'erie pr\'ec\'edente vue dans
$\Frac(\CC/\QQ)[[x]]$.

\begin{theoI}\label{t}
Soit $b(s) \in \Frac(\CC/\QQ)[s]$ le polyn\^ome de Bernstein (formel)
de $(f)_\QQ$ alors~:
\begin{description}
\item[(i)]
$b(s)$ est \`a racines rationnelles.
\item[(ii)]
$b(s)$ est le polyn\^ome unitaire de plus bas degr\'e dans $\C[s]$ tel
qu'il existe $h(x,y) \in \O_Z$ tel que
\begin{equation}\label{eq:faible}
\begin{cases}
h(x,y) \cdot b(s) f^s \in \D_{Z/Y}[s] \cdot f^{s+1} + \QQ \cdot \D_{Z/Y}
\sdt \cdot f^s\\
\textrm{ avec } h(0,y) \in \O_Y \smallsetminus \QQ.
\end{cases}
\end{equation}
\item[(iii)]
Il existe $h' \in \O_Y \smallsetminus \QQ$ tel que pour tout $y \in
V(\QQ) \smallsetminus V(h')$, $b(s)$ est le polyn\^ome de Bernstein
local (en $x=0$) de $f_y$.
\end{description}
\end{theoI}

\begin{remI}\label{rem:theo}
\begin{itemize}
\item[(a)]
Le polyn\^ome de Bernstein $b_g$ d'un $g\in \k[[x]]$ ($\k$ \'etant un
corps de caract\'eristique nulle) est non nul (Bj\"ork \cite{bjorkP})
ainsi notre polyn\^ome $b(s)$ est non nul. Par contre le fait que
$b_g$ soit \`a racines rationnelles est \`a notre connaissance une
question ouverte donc le point {\bf (i)} du th\'eor\`eme n\'ecessite
une d\'emonstration.
\item[(b)]
Dans {\bf (ii)}, la sp\'ecialisation en $y_0 \in V(\QQ)\smallsetminus
V(h(0,y))$ fait de $b(s)$ \emph{un} polyn\^ome de Bernstein de
$f_{y_0}$. D'apr\`es {\bf (iii)}, c'est en fait \emph{le} polyn\^ome
de Bernstein de $f_{y_0}$ si de plus $h'(y_0)\ne 0$.
\item[(c)]
La relation (\ref{eq:faible}) est en g\'en\'eral fausse si l'on
remplace $\D_{Z/Y} \sdt$ par $\D_{Z/Y}[s]$ (voir la d\'emonstration
en section 2).
\item[(d)]
Supposons $\QQ=(0)$. Le point {\bf (ii)} \'etablit le fait que
le polyn\^ome $b(s)$ (polyn\^ome de Bernstein formel de $f$ vue dans
$\Frac(\O_Y)[[x]]$) est solution de l'\'equation (\ref{eq:faible}) (et
c'est en fait dans $\C[s]$ le g\'en\'erateur de l'id\'eal des
solutions).
Dans \cite[Sect. 2.7]{bioscaT}, H. Biosca se pose le probl\`eme
inverse. Elle se donne l'\'equation (\ref{eq:faible}) et en cherche
une solution (en fait elle travaille avec plusieurs fonctions $f_j$)
et elle montre qu'un it\'er\'e d'un polyn\^ome de Bernstein absolu
local en $(x,y)=0$ est une solution.
En cons\'equence de {\bf (ii)}, nous obtenons~: il existe $N\in \N$
tel que le polyn\^ome $b(s)$ divise $b_f(s) \cdots b_f(s+N)$ ce qui au
passage implique le point {\bf (i)} de notre th\'eor\`eme (gr\^ace \`a
\cite{kashiwara}). Ici $b_f$ est le polyn\^ome de Bernstein (absolu)
de $f$.
\item[(e)]
D'un point de vue g\'eom\'etrique, $b(s)$ n'est int\'eressant que
si $f$ n'est pas g\'en\'eriquement lisse sur $V(\QQ)$ sinon d'apr\`es
{\bf (iii)} :
$b(s)=s+1$ si $f_y(0)=0$ pour $y$ g\'en\'erique dans $V(\QQ)$ et
$b(s)=1$ sinon.
\end{itemize}
\end{remI}

Si $f$ est polynomiale, nous pouvons pr\'eciser le point {\bf (ii)} du
th\'eor\`eme pr\'ec\'edent. De fa\c{c}on plus g\'en\'erale, soit $\CC$ un
anneau commutatif, int\`egre, unitaire et contenant les
nombres rationnels et soit $f\in \CC[x]$ (c'est le cadre de \cite{PJA}).
Notons $\An(\CC)$ l'alg\`ebre de Weyl au dessus de $\CC$ et posons
$\An(\CC)\sdt=\An(\CC) \otimes_\Q \Q\sdt$. Pour $\QQ \in \spec(\CC)$,
consid\'erons $(f)_\QQ$ qui est donc dans $\Frac(\CC/\QQ)[x]$.
Puisque $\Frac(\CC/\QQ)$ est de caract\'eristique nulle (car $\Q$ est
inclus dans $\CC$), $b(s)$ le polyn\^ome de Bernstein formel de $(f)_\QQ$
est \`a racines rationnelles (voir \cite{briancon}, voir aussi \cite{PJA}).

\begin{propI}\label{p}
Le polyn\^ome $b(s)$ est le polyn\^ome unitaire de plus bas degr\'e dans
$\C[s]$ tel que
\begin{equation}\label{eq:polyL}
\begin{cases}
h(x) \cdot b(s) f^s \in \An(\CC)[s] \cdot f^{s+1} +
\QQ \cdot \An(\CC)\sdt \cdot f^s\\
\textrm{ avec } h(x) \in \CC[x] \textrm{ et } h(0)\in \CC \smallsetminus \QQ.
\end{cases}
\end{equation}
\end{propI}

\begin{remI}\label{r}
Notons $b_{glob}(s)$ le polyn\^ome de Bernstein global de $(f)_\QQ$ alors
c'est le polyn\^ome unitaire de plus bas degr\'e dans $\C[s]$ r\'ealisant
\begin{equation}\label{eq:polyG}
\begin{cases}
h \cdot b_{glob}(s) f^s \in \An(\CC)[s] \cdot f^{s+1} +
\QQ \cdot \An(\CC)\sdt \cdot f^s\\
\textrm{ avec } h \in \CC \smallsetminus \QQ.
\end{cases}
\end{equation}
Autrement dit, c'est \emph{le} polyn\^ome de Bernstein g\'en\'erique global
de $f$ sur $V(\QQ)$ ce qui compl\`ete \cite{PJA}.
\end{remI}

Revenons \`a la situation de d\'epart~: $f\in \O_{X \times Y}$.
Comme cons\'equence du th\'eor\`eme \ref{t}, nous obtenons le r\'esultat
de constructibilit\'e suivant~:

\begin{corI}\label{c1}
La partition de $Y$ d\'efinie par le polyn\^ome de Bernstein de $f_y$
est constructible.
\end{corI}

Voici le dernier des principaux r\'esultats que nous d\'emontrerons.

\begin{corI}\label{c2}
Supposons que g\'en\'eriquement sur $V(\QQ)$, le lieu singulier relatif
de $f$ se projette sur $0$ par la projection canonique $Z\to X$.
Alors le polyn\^ome $b(s)$ du th\'eor\`eme \ref{t} est le polyn\^ome
unitaire de plus bas degr\'e dans $\C[s]$ pour lequel il existe $H \in \O_Y
\smallsetminus \QQ$ tel que:
\begin{equation}\label{eq:fort}
H(y) \cdot b(s) f^s \in \D_{Z/Y}[s] \cdot f^{s+1} + \QQ \cdot \D_{Z/Y}
\sdt \cdot f^s.
\end{equation}
\end{corI}
Si $\QQ=(0)$, cela nous dit que $b(s)$ est \emph{le} polyn\^ome de
Bernstein g\'en\'erique de $f$ au sens de Brian\c{c}on et al. \cite{bgm}.
Remarquons cependant que dans ce corollaire, on ne suppose rien sur $f_0$.
On est donc en dehors du cadre \'etudi\'e dans loc. cit.

Pour finir, d\'ecrivons la structure de l'article. Dans un premier temps
(section 2), nous d\'emontrons la proposition \ref{p} et la remarque \ref{r}.
Nous d\'emontrons aussi les corollaires \ref{c1} et \ref{c2} et
l'assertion (c) de la remarque \ref{rem:theo} supposant aquis le
th\'eor\`eme \ref{t}.
Le reste du papier est consacr\'e \`a le d\'emontrer. Nous rappelons
(section 3) ce qui nous sera n\'ecessaire concernant les bases standard
g\'en\'eriques \cite{genSB}, puis nous donnons (section 4) un algorithme
(infini) de calcul du polyn\^ome de Bernstein formel pour $f\in \k[[x]]$
(il s'agit d'une g\'en\'eralisation d'un algorithme de T.~Oaku \cite{oaku}).
En section 5, nous d\'ebutons la preuve du Th. \ref{t} qui consiste \`a
suivre l'algorithme pas \`a pas. Une premi\`ere \'etape consiste en
l'\'elimination de variables globales que sont $\ddt$ et les $\dx{i}$,
la seconde, plus technique, est une ``\'elimination'' des variables locales
$x_i$. \`A la fin de la premi\`ere \'etape nous serons en mesure de
d\'emontrer {\bf (i)} puis {\bf (iii)}. La derni\`ere section est
consacr\'ee \`a la preuve de {\bf (ii)}; cela passe par l'\'etablissement
d'une relation fonctionnelle g\'en\'erique formelle puis un ``passage du
formel \`a l'analytique'' inspir\'e de \cite{bmai90}.

\section{Premi\`eres d\'emonstrations}

Ici nous d\'emontrons des r\'esultats dont la preuve est ind\'ependante
du reste du papier: la proposition \ref{p}, la remarque \ref{r}; et
supposant acquis le th\'eor\`eme \ref{t}, les deux
corollaires ainsi que le (c) de la remarque \ref{rem:theo}.

\begin{proof}[D\'emonstration de la proposition \ref{p}]
Si $c(s)\in \C[s]$ satisfait (\ref{eq:polyL}) alors il est un polyn\^ome
de Bernstein local de $(f)_\QQ$ et est donc multiple de $b(s)$ donc
il suffit de montrer que $b(s)$ satisfait (\ref{eq:polyL}).
Il est bien connu que le polyn\^ome de Bernstein formel d'un polyn\^ome
$g\in \k[x]$ est le plus petit polyn\^ome unitaire tel que $q(x) b(s) g^s
\in \An(\k)[s] g^{s+1}$ avec $q(x) \in \k[x]$ et $q(0)\ne 0$ (voir par
exemple \cite{bmai90}).
Appliquons ceci \`a $(f)_\QQ \in \FF(\QQ)[x]$. En cons\'equence, il existe
$c\in \CC \smallsetminus \QQ$, $q' \in \CC[x]$ avec $q(0) \notin \QQ$ et
$P\in \An(\CC)[s]$ tels que~: $\big( c \cdot q'(x) \cdot b(s) -P f
\big)_\QQ \cdot ((f)_\QQ)^s =0$.

Posons $U=c \cdot q'(x) \cdot b(s) -P f \in \An(\CC)[s]$ et notons
$I'$ l'id\'eal de $\An(\CC)\sdt$ engendr\'e par $s + f(x)\ddt$ et $\dx{i}
+ \dxsur{f}{i}\ddt$ pour $i=1,\ldots,n$.
Son sp\'ecialis\'e $(I')_\QQ \subset \An(\FF(\QQ))\sdt$ est l'annulateur de
$(f)_\QQ^s$ (voir \cite{bmai02}). Ainsi $(U)_\QQ \in (I')_\QQ$.
\'Ecrivons $(U)_\QQ=\sum_j u_j \cdot (m_j)_\QQ$ avec $m_j\in I'$. Dans
cette \'ecriture, on rel\`eve les $u_j$, on chasse les d\'enominateurs
et on obtient l'existence de $c' \in \CC \smallsetminus \QQ$ tel que
\[c c' \cdot q'(x) \cdot b(s) \in \An(\CC)[s] f + I' +\An(\QQ)\sdt.\]
En appliquant cet op\'erateur \`a $f^s$, nous obtenons la relation voulue
ce qui d\'emontre la proposition.
\end{proof}

Pour la remarque \ref{r}, les arguments sont tout \`a fait similaires.
Nous laissons les d\'etails au lecteur. Maintenant supposons acquis
le th\'eor\`eme \ref{t} et commen\c{c}ons par une

\begin{proof}[Esquisse de preuve du Cor. \ref{c1}]
On montre que pour tout ferm\'e de Zariski $V$ de $Y$, on a une stratification
$V=\cup W$ en espaces localement ferm\'es telle que sur chaque strate le
polyn\^ome de Bernstein de $f_y$ est constant. Cela se fait par r\'ecurrence
sur la dimension de $V$ (le r\'esultat \'etant trivial si $\dim V=0$).
On \'ecrit $V=V(\QQ_1) \cup \cdots \cup V(\QQ_p)$ comme r\'eunion de ses
composantes irr\'eductibles et on applique le {\bf (iii)} du Th.~\ref{t}
\`a chacun des id\'eaux premiers $\QQ_i$, on note $h_i$ le $h'$ obtenu.
On a alors $V=V_1 \cup V_2$ o\`u $V_1=\cup (V(\QQ_i) \smallsetminus V(h_i))$
et sur chaque strate de $V_1$ le polyn\^ome de Bernstein est constant.
On applique l'hypoth\`ese de r\'ecurrence \`a $V_2=\cup (V(\QQ_i) \cap 
V(h_i))$ qui est de dimension strictement inf\'erieure \`a $\dim V$.
\end{proof}
\noindent
Maintenant, d\'emontrons le corollaire \ref{c2}.

\begin{proof}
Cette preuve s'inspire de celle de \cite[Prop. 1.4]{bgm}.
Pour d\'emontrer le corollaire, il suffit de montrer que $b(s)$
satisfait (\ref{eq:fort}). Notons $J$ l'id\'eal de $\O_Z$ engendr\'e par
$f$ et les $\frac{\partial f}{\partial x_i}$. Par hypoth\`ese, il existe
$h_0 \in \O_Y \smallsetminus \QQ$ tel que pour tout $y_0 \in V(\QQ)
\smallsetminus V(h_0)$, $V(J_{|y_0})= \{0\} \subset X$.
Par cons\'equent, $V(J+\O_Z \cdot \QQ) \smallsetminus
V(\O_Z \cdot h_0) \subset \{0\} \times V(\QQ)$. Autrement dit
\[V(\sqrt{J+\O_Z \cdot \QQ} : h_0) \subset \{0\} \times V(\QQ).\]
Ainsi le lieu des z\'eros de $(\sqrt{J+\O_Z \cdot \QQ} : h_0)$ et de $h$
(celui de {\bf (ii)} dans le Th.~\ref{t}) est inclus dans le lieu des
z\'eros de $h(0,y)$. Ainsi, pour $l\in \N$ assez grand, $h_1:=h(0,y)^l$ est
dans l'id\'eal  $(\sqrt{J+\O_Z \cdot \QQ} : h_0) + \O_Z \cdot h$.
On constate alors que pour un certain $k\in \N$, $H:=(h_0 h_1)^k$ appartient
\`a l'id\'eal $J+ \O_Z \cdot \QQ +\O_Z \cdot h$. Notons que $H \in \O_Z
\smallsetminus \QQ$. Maintenant, puisque $f_y(0)=0$ pour $y$ g\'en\'erique
dans $V(\QQ)$ alors $b(s)$ est multiple de $(s+1)$. Ainsi, si l'on note
$b(s)=(s+1) \tilde{b}(s)$, on a pour tout $i$,
\[b(s) \frac{\partial f}{\partial x_i} f^s = \tilde{b}(s)
\frac{\partial}{\partial x_i} \cdot f^{s+1}.\]
De cette \'equation et de la relation (\ref{eq:faible}), on obtient la
relation d\'esir\'ee (\ref{eq:fort}).
\end{proof}

Pour finir, d\'emontrons l'assertion (c) de la remarque \ref{rem:theo}.

\begin{proof}
Soit $f$ une fonction analytique de $n+1$ variables complexes
$x_1,\ldots,x_n,y$. Supposons que l'hypoth\`ese du corollaire \ref{c2}
soit v\'erifi\'ee pour $\QQ=(0)$. De plus supposons que le nombre de Milnor
de $f_0$ soit diff\'erent de celui de $f_{y_0}$ pour $y_0\ne 0$ proche de
$0$. Par exemple, on peut prendre $f=x_1^2+y x_2^2 +x_2^3$. Sur cet
exemple, on a $V(f, \frac{\partial f}{\partial x_1},
\frac{\partial f}{\partial x_2})=\{ (0,0)\} \times \C$.
Appliquons le corollaire \ref{c2} \`a $\QQ=(0)$~:
\[(\star) \qquad H(y) \cdot b(s) f^s \in \D_{Z/Y}[s] \cdot f^{s+1}.\]
Quitte \`a multiplier cette relation par une unit\'e de $\O_{Y,0}$ on peut
supposer que $H(y)=y^N$ pour un certain entier $N$. Maintenant appliquons
l'assertion {\bf (ii)} du th\'eor\`eme \ref{t} \`a $\QQ=(y)$. Notons $b_0(s)$
le polyn\^ome de Bernstein en question et (par l'absurde) supposons fausse
l'assertion (c) de la remarque \ref{rem:theo}, on obtient donc~:
\[(\star\star) \qquad h(x,y) \cdot b_0(s) f^s \in \D_{Z/Y}[s] \cdot f^{s+1}+
y \cdot \D_{Z/Y}[s] \cdot f^s\]
avec $h(0,y)$ ne s'annulant pas en $0$. Cette derni\`ere condition nous
dit que $h$ est inversible, on peut donc le supposer \'egal \`a $1$.
En it\'erant $(\star\star)$ on obtient
$b_0(s)^N f^s \in \D_{Z/Y}[s]\cdot f^{s+1} + y^N\cdot \D_{Z/Y}[s]\cdot f^s$.
En multipliant cette derni\`ere relation par $b(s)$ et en utilisant $(\star)$,
on obtient $b(s) b_0(s)^N f^s \in  \D_{Z/Y}[s]\cdot f^{s+1}$ i.e. il existe
un polyn\^ome de Bernstein relatif non nul. Or $f$ est une d\'eformation
\`a un param\`etre de $f(x,0)$ \`a nombre de Milnor non constant ce qui
contredit \cite[Th. 4]{blm}.
\end{proof}


\section{Bases standard param\'etriques}

\subsection{Bases standard g\'en\'eriques}

Pour que le papier soit le plus autonome possible, nous avons d\'ecid\'e
de donner les rappels n\'ecessaires concernant les bases standard.
Cependant, afin de garder une taille raisonable \`a ce papier, nous
n'entrerons pas dans tous les d\'etails. Le lecteur sera renvoy\'e \`a
(Castro-Jim\'enez, Granger \cite{cg}).
Nous donnons ensuite les r\'esultats n\'ecessaires sur les
bases standard g\'en\'eriques (voir \cite{genSB} pour un traitement
plus complet). Enfin dans le paragraphe suivant, nous appliquons ces
derni\`eres \`a l'\'elimination g\'en\'erique de variables globales.

Dans la suite, nous aurons besoin de travailler dans plusieurs types
de $\k[[x]]$-alg\`ebres non commutatives ($\k$ \'etant un corps de
caract\'eristique $0$). Dans cette section, nous donnons une construction
qui couvre tous les cas rencontr\'es plus loin. Ce qui suit peut \^etre vu
comme une version locale proche de \cite{bmai02}.

Soit $x=(x_1,\ldots,x_n)$ et $z=(z_1,\ldots, z_q)$ deux syst\`emes de
variables. Soit $\ring=\ring(\k):=\k[[x]]\langle z \rangle$ la
$\k[[x]]$-alg\`ebre engendr\'ee par les $z_i$ avec les relations
de commutation suivantes~:
\begin{itemize}
\item[(i)]
$\dps [z_i, a(x)]\in \k[[x]]$ pour $a(x)\in \k[[x]]$,
\item[(ii)]
$\dps [z_i,z_j] \in \k[[x]]+ \sum_{k=1}^q \k[[x]] z_k$.
\end{itemize}

La notation $\k[[x]]\langle z \rangle$ rappelle que les variables $z_i$
ne sont pas commutatives en g\'en\'eral.
Les cas que nous rencontrerons dans la suite sont: $\FDn(\k)\sdt=
\FDn(\k) \otimes_\k \k\sdt$ o\`u $\FDn(\k)$ est l'anneau des op\'erateurs
diff\'erentiels \`a coefficients dans $\k[[x]]$, $\FDn(\k)[s]$, $\k[[x]][s]$,
$\k[s]$. Tous ces cas sont couverts par la construction ci-dessus.

\begin{rem}\label{rem:uniq_ecriture}
Avec cette d\'efinition de $\ring$, on n'a pas n\'ecessairement
unicit\'e de l'\'ecriture \`a gauche. Par exemple, si $\ring=\k
\langle z_1, z_2,z_3 \rangle$ ($n=0$, $q=3$),  $[z_1, z_2]=z_1$,
$[z_1, z_3]=z_2$ et $[z_2, z_3]=z_1$ alors $z_3 z_1 z_2$ aura (au
moins) deux \'ecritures possibles~: $z_1 z_2 z_3- z_1^2 -z_2^2$ et $z_1
z_2 z_3- z_1^2 -z_2^2 -z_2$. On obtient la premi\`ere en faisant
commuter (dans le terme de degr\'e $3$) deux $z_i$ suivant les
transpositions d'indices~: $(3,1)$ et $(3,2)$. Pour la seconde, on utilise
$(1,2)$ suivie de $(3,2)$, $(3,1)$ et $(2,1)$.
\end{rem}

Dans la suite, nous imposons donc l'hypoth\`ese suppl\'ementaire
d'unicit\'e de l'\'ecriture \`a gauche (cette unicit\'e, bien entendu,
a lieu dans tous les anneaux \'enum\'er\'es ci-dessus).

Pour commencer, nous devons \'enoncer un th\'eor\`eme de division dans
$\ring$.

Pour $P\in \ring$ s'\'ecrivant (de mani\`ere unique) $P=\sum c_{\alpha
  \beta} x^\alpha z^\beta$, $(\alpha, \beta) \in \N^{n+q}$, $c_{\alpha
  \beta}\in \k$, on d\'efinit son diagramme de Newton $\DN(P) \subset
\N^{n+q}$ comme l'ensemble des $(\alpha, \beta)$ tels que $c_{\alpha
  \beta}$ est non nul.

Soit $\prec$ un ordre (total et compatible avec l'addition) sur les
$(\alpha, \beta) \in \N^{n+q}$ d\'efini comme suit~:
on se donne une forme lin\'eaire $L(\beta)=\sum_i l_i \beta_i$, tels que
les $l_i$ soient positifs ou nuls et on d\'efinit $\prec=\prec_L$~:
\[(\alpha, \beta) \prec (\alpha', \beta') \iff
\begin{cases}
L(\beta) < L(\beta')\\
\textrm{ou \'egalit\'e et } |\beta|< |\beta'|\\
\textrm{ou \'egalit\'es et } |\alpha| > |\alpha'|\\
\textrm{ou \'egalit\'es et } (\alpha, \beta) >_0 (\alpha', \beta').
\end{cases}\]
Ici $<_0$ est un ordre total, bon et compatible avec l'addition dans
$\N^{n+q}$.

Pour $P\in \ring$ non nul, on note $\exp_\prec(P)$ le maximum de
$\DN(P)$ pour $\prec$. C'est son exposant privil\'egi\'e. On note aussi son
terme et coefficient privil\'egi\'e~: $\tp_\prec(P)=(x, z)^{\exp_\prec(P)}$,
$\cp_\prec(P)=c_{\exp_\prec(P)}$.

Soient $P_1,\ldots,P_r \in \ring$. On d\'efinit une partition
$\N^{n+q}=\Delta_1\ \cup \cdots\ \cup \Delta_r \cup \bar{\Delta}$
associ\'ee aux $\exp_\prec(P_j)$ comme suit: $\Delta_1= \exp_\prec(P_1)+
\N^{n+q}$, puis pour $j\ge 2$, $\Delta_j=(\exp_\prec(P_j) +\N^{n+q})
\smallsetminus \cup_{k=1}^{j-1} \Delta_k$.

\begin{theo}[Th\'eor\`eme de division]\label{theo:div}
Pour $P \in \ring$, il existe un unique $(Q_1,\ldots,Q_r,R)\in \ring^{r+1}$
tel que
\begin{itemize}
\item
$P=\sum_j Q_j P_j +R$
\item
pour tout $j$, $Q_j=0$ ou bien $\DN(Q_j)+\exp_\prec(P_j) \subset \Delta_j$,
\item
$R=0$ ou bien $\DN(R) \subset \bar{\Delta}$.
\end{itemize}
\end{theo}

\begin{proof}[Id\'ee de la preuve]
L'unicit\'e est facile, occupons nous de l'existence. Pour cela, nous allons
nous ramener aux r\'esultats de \cite{cg}.
\begin{lem}
Il existe une forme lin\'eaire $L'$ \`a coefficients strictement positifs
agissant sur les $\beta$ tel que pour $j=1,\ldots,r$,
$\exp_{\prec_L} (P_j)=\exp_{\prec_{L'}}(P_j)$
\end{lem}
La preuve de ce lemme se fait exactement comme celle de
\cite[Prop. 8]{acg01}.
Remarquons que la division ne d\'epend que des $\exp_\prec
(P_j)$, ainsi gr\^ace \`a ce lemme nous pouvons supposer que $L$ est \`a
coefficients strictement positifs.
La forme $L$ donne lieu \`a une filtration sur $\ring$ dont le gradu\'e
est isomorphe \`a $\k[[x]][\xi_1,\ldots,\xi_q]$. Ici les $\xi_i$ sont
des variables commutatives correspondant aux $z_i$. En consid\'erant
les symboles principaux de $P$ et des $P_j$ par rapport \`a $L$, on se
ram\`ene \`a une division dans cet anneau. Ainsi,
la preuve se fait exactement comme celle de \cite[Th. 2.4.1]{cg}.
\end{proof}

\noindent
Voici quelques d\'efinitions et r\'esultats utiles pour la suite (voir
\cite{cg} pour les d\'emonstrations).
\begin{itemize}
\item[1.]\label{rappel1}
Dans le th\'eor\`eme pr\'ec\'edent, on a :

$\exp_{\prec_L}(P)=\max \{ \exp_{\prec_L}(Q_j P_j), j=1,\ldots, r;
\exp_{\prec_L}(R)\}$.\\
On d\'efinit $\ord^L(P)$ comme \'etant le maximum des $L(\beta)$
pour $(\alpha, \beta)\in \DN(P)$.
En cons\'equence,

$\ord^L(P)=\max\{\ord^L(Q_j P_j), j=1,\ldots, r; \ord^L(R)\}$.
\item[2.]\label{rappel2}
Pour un id\'eal $J \subset \ring$, on d\'efinit $\Exp_\prec(J)$
comme l'ensemble des $\exp_\prec(P)$ pour $P\in J$ non nul. Cet ensemble
est stable par addition dans $\N^{n+q}$, ainsi il existe $P_1,\ldots,P_r$
dans $J$ tels que $\Exp_\prec(J)=\bigcup_j (\exp_\prec(P_j) +\N^{n+q})$.
Un tel ensemble est appel\'e \emph{base standard de} J (pour $\prec$).
\item[3.]\label{rappel3}
Soient $P_1,\ldots, P_r \in J$. Les assertions suivantes sont
\'equivalentes.
\begin{itemize}
\item
$P_1,\ldots,P_r$ forment une $\prec$-base standard de $J$.
\item
Pour $P\in \ring$~: $P\in J \iff$ le reste $R$ de la division de $P$
par les $P_j$ est nul.
\end{itemize}
\item[4.]\label{rappel4}
$S$-op\'erateurs et crit\`ere de Buchberger (\cite{buchberger}
dans le cas polynomial).
\begin{itemize}
\item
Soient $P,P'\in \ring$. Notons $e=\exp_\prec(P)$ et $e'=
\exp_\prec(P')$. Soit $\mu=\max(e,e')$ que l'on d\'efinit en posant $\mu_i=
\max(e_i,e'_i)$ pour chaque $i=1,\ldots,n+q$. On d\'efinit alors le
$S$-op\'erateur de $P$ et $P'$~: $S(P,P')= \cp_\prec(P') m P-
\cp_\prec(P)m' P'$ o\`u $m=(x,z)^{\mu-e}$ et $m'=(x,z)^{\mu-e'}$.
\item
Soit $\GG$ un syst\`eme de g\'en\'erateurs de $J\subset \ring$, alors
$\GG$ est une base standard de $J$ si pour tout $P,P'\in \GG$, le reste de
la division de $S(P,P')$ par $\GG$ est nul \cite[Prop. 2.5.1]{cg}.
\end{itemize}
\end{itemize}
\

Maintenant introduisons les bases standard g\'en\'eriques. Soit $\CC$ un
anneau int\`egre commutatif unitaire et contenant comme sous-anneau le corps
des nombres rationnels. En ce qui nous concerne, il faut penser \`a
$\CC=\O_Y$, cependant $\CC$ peut \^etre \'egal \`a d'autres anneaux, tel
que $\k[y]$. Soit $\FF=\Frac(\CC)$ son corps des fractions.
On note $\spec(\CC)$ et $\specm(\CC)$ son spectre et son spectre maximal,
respectivement. Dans la suite, lorsque $\CC=\O_Y$, nous identifierons
$\specm(\CC)$ et $Y$.
Pour tout $\PP \in \spec(\CC)$, $\FF(\PP)$ d\'esigne le corps des fractions
de $\CC/\PP$; c'est un corps de caract\'eristique $0$ (ceci gr\^ace \`a
l'hypoth\`ese $\Q \subset \CC$). Pour tout id\'eal $\II$, $V(\II)=\{\PP
\in \spec(\CC)| \II \subset \PP\}$ d\'esigne le ferm\'e de Zariski d\'efini
par $\II$. On notera $V_m(\II)$ sa restriction \`a $\specm(\CC)$.

Soit $\PP \in \spec(\CC)$ et $c\in \CC$. On note $[c]_\PP$ sa classe dans
$\CC/\PP$ et $(c)_\PP=\frac{[c]_\PP}{1}$ cette classe vue dans le corps
$\FF(\PP)$. On appelle $(\cdot)_\PP$ la sp\'ecialisation en $\PP$.\\

Dans ce qui suit, afin de rendre l'exposition plus rigoureuse, nous
invoquons le language des cat\'egories mais vue la simplicit\'e de
notre situation, nous aurions pu l'\'eviter.
Consid\'erons la cat\'egorie dont un objet est $A[[x]]$ o\`u $A$ est
un anneau et les fl\`eches sont des applications (ensemblistes).

On se donne une fl\`eche $\phi$ de l'objet $\CC[[x]]$ vers lui-m\^eme
et pour tout $\PP \in \spec(\CC)$, $\phi_\PP$ une fl\`eche de
$\FF(\PP)[[x]]$ vers lui m\^eme. Nous dirons que $\phi$ est adapt\'ee
aux $\phi_\PP$ si pour tout $a(x)\in \CC[[x]]$, $\big( \phi(a(x))
\big)_\PP=\phi_\PP \big( (a(x))_\PP \big)$.


Maintenant, pour tout $\PP$, on se donne $\ring(\FF(\PP))=\FF(\PP)[[x]]
\langle z \rangle$. Pour chaque $i=1,\ldots,q$, on d\'efinit
la fl\`eche $\phi_{i,\PP}$ de $\FF(\PP)[[x]]$ vers lui-m\^eme en posant
$\phi_{i,\PP}(c(x)):=[z_i, c(x)]$.
On d\'efinit alors $\ring(\CC)$ comme la $\CC$-alg\`ebre
engendr\'ee par $\CC[[x]]$ et $z_1,\ldots,z_q$ avec les relations
de commutation suivantes:
\begin{itemize}
\item[(i)]
$\dps [z_i, a(x)]=\phi_i(a(x))$ pour $a(x)\in \CC[[x]]$,
\item[(ii)]
$\dps [z_i,z_j] = u_{ij}+ \sum_{k=1}^q v_{ijk} z_k$.
\end{itemize}
Ici on a fix\'e les $u_{ij}$ et $v_{ijk}$ comme \'etant des entiers
dans $\Z$ (pour simplifier) et $\phi_i$ est une fl\`eche adapt\'ee
aux $\phi_{i,\PP}$. De plus on suppose que les commutateurs $[z_i,
z_j]$ soient les m\^emes dans $\ring(\CC)$ et dans les
$\ring(\FF(\PP))$.
Dans la suite, les choses seront simples et l'on aura
$\phi_i=\frac{\partial}{\partial x_i}$ lorsque $z_i$ sera une
d\'erivation partielle et $\phi_i=0$ dans les autres situations.

Une fois cette d\'efinition faite, on \'etend de fa\c{c}on naturelle 
les op\'erations de sp\'ecialisation aux \'el\'ements de $\ring(\CC)$,
ainsi qu'\`a ceux de $\ring(\FF)$, dont le d\'enominateur des
coefficients n'est pas dans $\PP$.

Maintenant si $J$ est un id\'eal de $\ring(\CC)$, on note $(J)_\PP$
l'id\'eal de $\ring(\FF(\PP))$ engendr\'e par les $(P)_\PP$ avec $P\in J$.

Fixons un ordre $\prec=\prec_L$ comme plus haut et fixons $\QQ \in
\spec(\CC)$. On note $\ring(\QQ)$ l'id\'eal de $\ring(\CC)$ dont les
\'el\'ements ont leurs coefficients dans $\QQ$.

Soit $P\in \ring(\CC) \smallsetminus \ring(\QQ)$, qu'on \'ecrit comme
plus haut sauf qu'ici les $c_{\alpha \beta}$ sont dans $\CC$.
On d\'efinit $\DN^\modQ(P)$ son diagramme de Newton modulo $\QQ$ comme
l'ensemble des $(\alpha, \beta)$ tels que $c_{\alpha \beta}\in \CC
\smallsetminus \QQ$.
En fait $\DN^\modQ(P)=\DN((P)_\QQ)$. On d\'efinit $\exp_\prec^\modQ(P)$ son
exposant privil\'egi\'e modulo $\QQ$ comme le maximum (pour $\prec$) de
$\DN^\modQ(P)$. On d\'efinit aussi son terme et son coefficient
privil\'egi\'e modulo $\QQ$~: $\tp_\prec^\modQ(P)=
(x,z)^{\exp_\prec^\modQ(P)}$, $\cp_\prec^\modQ(P)=c_{\exp_\prec^\modQ(P)}$.

Soit $J$ un id\'eal de $\ring(\CC)$. Soit $\Exp_\prec^\modQ(J)$
l'ensemble des $\exp_\prec^\modQ(P)$ pour $P \in J \smallsetminus
\ring(\QQ)$. Il est facile de voir que cet ensemble est stable par
addition. Ainsi, par le lemme de Dickson, la d\'efinition suivante n'est
pas vide.

\begin{defin}
On d\'efinit une base standard g\'en\'erique de $J$ sur $V(\QQ)$ (pour
$\prec$) comme un sous ensemble fini $\GG=\{P_1,\ldots,P_r\}$ de $J$ tel que
$\Exp_\prec^\modQ(J) =\bigcup_j (\exp_\prec^\modQ(P_j)+ \N^{n+q})$.
\end{defin}
Dans \cite{genSB}, nous avons donn\'e une d\'efinition plus g\'en\'erale.
Cependant pour l'usage qu'on en fera ici, la d\'efinition ci-dessus est
suffisante.

Notons $\langle \QQ \rangle$ l'id\'eal de $\ring(\FF)$ constitu\'e
d'\'el\'ements dont le num\'erateur des coefficients est dans $\QQ$.

\begin{prop}[Division modulo $\QQ$, {\cite[Prop. 2.1.2]{genSB}}]
Soient $P_1,\ldots,P_r \in \ring(\CC)$ et soit $\Delta_1\cup \cdots
\cup \Delta_r \cup \bar{\Delta}$ la partition de $\N^{n+q}$ associ\'ee aux
$\exp_\prec^\modQ(P_j)$. Pour $P\in \ring(\CC)$, il existe $Q_1,\ldots,
Q_r, R\in \ring(\FF)$ et $T \in \langle \QQ \rangle$ tels que
$P=\sum_j Q_j P_j +R +T$ et
\begin{itemize}
\item
$\DN(Q_j) +\exp_\prec^\modQ(P_j) \subset \Delta_j$ si $Q_j \ne 0$,
\item
$\DN(R) \subset \bar{\Delta}$ si $R\ne 0$,
\item
le d\'enominateur des coefficients de $R$, $T$ et des $Q_j$ sont des
puissances de $h=\prod_j \cp_\prec^\modQ(P_j)$. Autrement dit, la division
a lieu dans $\ring(\CC[h^{-1}])$ (i.e. les coefficients sont dans le
localis\'e de $\CC$ par rapport \`a $h$).
\end{itemize}
De plus $(Q_1,\ldots,Q_r, R)$ est unique modulo $\langle \QQ \rangle$. On
appelle $R$ le reste modulo $\QQ$.
\end{prop}

\begin{proof}
La preuve consiste \`a poser $P_j=P_j^1 -P_j^2$ avec $P_j^2 \in \langle \QQ
\rangle$ et $\exp_\prec(P_j^1)=\exp_\prec^\modQ(P_j)$ et \`a effectuer
la division de $P$ par les $P_j^1$ dans $\ring(\FF)$. Voir les d\'etails
dans \cite[Prop. 2.1.2]{genSB}.
\end{proof}

\begin{cor}\label{cor:restemodQ}
Soit $P\in \ring(\CC)$ tel que $(P)_\QQ \in (J)_\QQ$ et soit $\GG$ une
$\prec$-base standard g\'en\'erique de $J$ sur $V(\QQ)$ alors le reste
modulo $\QQ$ de la division modulo $\QQ$ de $P$ par $\GG$ est nul.
\end{cor}
\begin{proof}
\'Ecrivons $P=\sum_j Q_j P_j +R+T$ comme dans la proposition.
Remarquons que $\exp_\prec^\modQ(P_j)=\exp_\prec((P_j)_\QQ)$ donc
la partition de $\N^{n+q}$ dans la proposition est \'egale \`a celle
associ\'ee aux $\exp_\prec((P_j)_\QQ)$. Sp\'ecialisons l'\'egalit\'e
pr\'ec\'edente en $\QQ$ (c'est possible puisque $h \notin \QQ$).
On obtient $(P)_\QQ=\sum_j (Q_j)_\QQ (P_j)_\QQ +(R)_\QQ$ avec
$\DN((Q_j)_\QQ)+\exp_\prec((P_j)_\QQ) \subset \Delta_j$ et 
$\DN((R)_\QQ)\subset \bar{\Delta}$. Ainsi l'\'egalit\'e pr\'ec\'edente
est le r\'esultat de la division de $(P)_\QQ$ par $(\GG)_\QQ$, or $(P)_\QQ
\in (J)_\QQ$ donc par le rappel (3.) page \pageref{rappel3}, $(R)_\QQ=0$
i.e. $R\in \langle \QQ \rangle$.
\end{proof}

\begin{theo}[\cite{genSB} Th.~2.1.6]\label{theo:BSgen}
Soit $\GG=\{P_1,\ldots, P_r\}$ une base standard g\'en\'erique de $J$ sur
$V(\QQ)$ et soit $h=\prod_j \cp_\prec^\modQ(P_j)$. Pour tout
$\PP\in V(\QQ)\smallsetminus V(h)$, $(\GG)_\PP$ est une base standard
de $(J)_\PP$.
\end{theo}
Remarquons que $\exp_\prec((P_j)_\PP)= \exp_\prec^\modQ(P_j)$ par d\'efinition
de $h$, par cons\'equent $\Exp_\prec((J)_\PP)$ est g\'en\'eriquement constant
et \'egal \`a $\Exp_\prec^\modQ(J)$. Remarquons aussi que $V(\QQ)
\smallsetminus V(h)$ n'est pas vide puisque $h\notin \QQ$.

\begin{proof}
Nous allons utiliser le crit\`ere de Buchberger (rappel (4.)
page~\pageref{rappel4}). Pour $P\in J$, effectuons sa division modulo $\QQ$
par $\GG$~: $P=\sum_j Q_j P_j +R +T$; le reste modulo $\QQ$ est nul par le
corollaire \ref{cor:restemodQ}. Comme dans la preuve de ce corollaire, nous
sp\'ecialisons cette division en $\PP \in V(\QQ)\smallsetminus V(h)$
(ce qui est possible car $h\notin \PP$). Ce que nous obtenons est la
division de $(P)_\PP$ par $(\GG)_\PP$, division dont le reste est nul
puisque $\QQ \subset \PP$. En cons\'equence $(\GG)_\PP$ engendre $(J)_\PP$.

Maintenant soient $P,P'$ dans $\GG$ et $S=\cp_\prec^\modQ(P') m P-
\cp_\prec^\modQ(P)m' P'$ o\`u l'on a pos\'e
$m=(x,z)^{\mu-\exp_\prec^\modQ(P)}$, $m'=(x,z)^{\mu-\exp_\prec^\modQ(P')}$
et $\mu=\max(\exp_\prec^\modQ(P), \exp_\prec^\modQ(P'))$. On constate que
pour tout $\PP \in V(\QQ) \smallsetminus V(h)$, $(S)_\PP=S((P)_\PP,
(P')_\PP)$. Par les m\^emes arguments que ci-dessus, on montre que la
division de $(S)_\PP$ par $(\GG)_\PP$ a un reste nul. On conclut \`a
l'aide du crit\`ere de Buchberger.
\end{proof}

\subsection{\'Elimination g\'en\'erique de variables globales}

Ce que nous appelons variables globales sont les $z_i$.
Par opposition les $x_i$ sont dites variables locales. Dans ce paragraphe,
nous montrons comment \'eliminer g\'en\'eriquement les variables
$z_{p+1}, \ldots, z_q$ avec $p<q$.

\'Enon\c{c}ons d'abord le r\'esultat dans le cas non param\'etrique.
Soit $J$ un id\'eal dans $\ring(\k)$. Soit $L$ la forme d\'efinie
par $L(\beta)=\sum_{p+1}^q \beta_i$. Grossi\`erement on met un
poids strictement positif (ici $1$) sur les variables \`a \'eliminer et
un poids nul sur les autres.
Notons $\prec=\prec_L$ l'ordre d\'efini par $L$ et soit
$\prec'$ sa restriction aux $(\alpha,\beta') \in \N^{n+p}$. Notons que cet
ordre est un ordre sur $\N^{n+p}$ associ\'e \`a la forme lin\'eaire
$L'(\beta_1,\ldots,\beta_p)=\sum_1^p \beta_i$ donc le th\'eor\`eme de
division dans $\k[[x]]\langle z' \rangle$, $z'=(z_1,\ldots,z_p)$, s'applique.

\begin{prop}
Soit $\GG$ une $\prec$-base standard de $J$ alors $\GG'=\GG \cap
\k[[x]]\langle z' \rangle$
est une $\prec'$-base standard de $J \cap \k[[x]]\langle z' \rangle$.
\end{prop}
On dit d'un tel ordre $\prec$ que c'est un \emph{ordre qui \'elimine}
les variables $z_{p+1}, \ldots,z_q$.
\begin{proof}
Les arguments sont standard et similaires \`a ceux de \cite[1.7]{cg}.
Soit $P\in J\cap \k[[x]]\langle z' \rangle$. Divisons $P$ par $\GG$ par
rapport \`a l'ordre $\prec$~: $P=\sum_j Q_j P_j$. Par hypoth\`ese sur $P$,
$\ord^L(P)=0$ donc (voir rappel (1.) page \pageref{rappel1})
pour tout $j$ tel que $Q_j\ne 0$, $\ord^L(Q_j)=0$ et
pour un tel $j$, $\ord^L(P_j)=0$ i.e. $Q_j$ et $P_j$ sont dans $\k[[x]]
\langle z' \rangle$.
On constate alors que la division pr\'ec\'edente est la division de $P$ par
$\GG'$ par rapport \`a $\prec'$, division pour laquelle le reste est nul.
On ach\`eve la preuve en utilisant le rappel (3.) page \pageref{rappel3}.
\end{proof}

Maintenant soit $J$ dans $\ring(\CC)$ et $\GG$ une $\prec$-base standard
g\'en\'erique de $J$ sur $V(\QQ)$.
\begin{prop}
Soit $\GG' \subset \GG$ l'ensemble des $P_j$ tels que
$\tp_\prec^\modQ(P_j)$ est dans $\CC[[x]]\langle z' \rangle$. Alors $\GG'$
est une $\prec'$-base standard g\'en\'erique de $J'=\big(J+ \ring(\QQ) \big)
\bigcap \CC[[x]]\langle z' \rangle$ sur $V(\QQ)$.
\end{prop}
\begin{proof}
Par d\'efinition d'une base standard g\'en\'erique, il suffit de montrer que
pour tout $P\in J'$, $\exp_{\prec'}((P)_\QQ)$ appartient \`a $\exp_{\prec'}
((P_j)_\QQ)+\N^{n+p}$ pour un certain $P_j \in \GG'$.
Pour un tel $P$, on a $(P)_\QQ \in (J)_\QQ \cap \FF(\QQ)[[x]]\langle z'
\rangle$.
Maintenant il est facile de voir que $(\GG')_\QQ$ est \'egal \`a $(\GG)_\QQ
\cap \FF(\QQ)[[x]]\langle z' \rangle$. Or par d\'efinition de $\GG$ est par
la proposition pr\'ec\'edente, ce  dernier ensemble est une $\prec'$-base
standard de $(J)_\QQ \cap \FF(\QQ)[[x]]\langle z' \rangle$ ce qui ach\`eve
notre d\'emonstration.
\end{proof}

\begin{cor}\label{cor:eliminegen}
Soit $h$ le produit des $\cp_\prec^\modQ(P_j)$ avec $P_j\in \GG$ alors pour
tout $\PP \in V(\QQ)\smallsetminus V(h)$,
\[(J)_\PP \cap \FF(\PP)[[x]]\langle z' \rangle =\Big(\big(J+ \ring(\QQ)
\big)\bigcap \CC[[x]]\langle z' \rangle \Big)_\PP\]
et ces id\'eaux sont engendr\'es par $(\GG')_\PP=(\GG)_\PP \cap
\FF(\PP)[[x]]\langle z' \rangle$.
\end{cor}
\begin{proof}
C'est une application directe du th\'eor\`eme \ref{theo:BSgen} et des deux
propositions pr\'ec\'edentes.
\end{proof}

\section{Construction algorithmique du polyn\^ome de Bernstein formel}

\'Etant donn\'ee une s\'erie formelle $f=f(x)\in \k[[x]]$ \`a $n$ variables
et \`a coefficients dans un corps $\k$ de caract\'eristique $0$,
J.~E.~Bj\"ork (\cite{bjorkP}, voir aussi \cite{bjorkB}) a d\'emontr\'e que
le polyn\^ome de Bernstein $b_f$ associ\'e est non nul. De plus si au
d\'epart $f\in \C\{x\}$ alors d'apr\`es J.~Brian\c{c}on et Ph.~Maisonobe
\cite{bmai90}, son polyn\^ome de Bernstein analytique est \'egal \`a son
polyn\^ome de Bernstein formel. Enfin toujours dans ce m\^eme cas,
M.~Kashiwara \cite{kashiwara} a d\'emontr\'e que les racines de $b_f$ sont
rationnelles n\'egatives. La rationalit\'e de $b_f$ pour $f\in \k[[x]]$
est, \`a notre connaissance, une question ouverte.

T.~Oaku \cite{oaku} a donn\'e un algorithme de calcul du polyn\^ome de
Bernstein formel pour $f\in \k[x]$. Cet algorithme se compose d'une premi\`ere
partie o\`u l'on \'elimine des variables globales et d'une seconde o\`u l'on
``\'elimine'' les variables locales $x_i$. Ici nous proposons une variante
de la premi\`ere partie (variante inspir\'ee de \cite{bmai02}) et montrons
que la seconde partie fonctionne pour $f\in \k[[x]]$.

Soit $f \in \k[[x]]$. Le module libre $\mathcal{L}=\k[[x]][1/f,s]\cdot f^s$
a une structure naturelle de $\FDn(\k)[s]$-module. Suivant B.~Malgrange
\cite{malgrange}, on en fait un $\FDnpun(\k)$-module (o\`u l'on batise $t$
la nouvelle variable)~: si $g(s)\in \k[[x]][1/f,s]$, on pose $t\cdot g(s)
f^s=g(s+1) f f^s$ et $\ddt \cdot g(s) f^s=-s g(s-1) f^{-1} f^s$. On constate
alors que $s$ agit sur $\mathcal{L}$ comme $-\ddt t$. Cette identification
permet de faire de $\mathcal{L}$ un $\FDn(\k)\sdt$-module (cette approche
est d\^ue \`a Brian\c{c}on et Maisonobe \cite{bmai02} dans le cas
alg\'ebrique).
Consid\'erons les id\'eaux suivants.

\begin{description}
\item[0]\label{debutalgo}
$\dps I(f)=\FDn(\k)\sdt\cdot (s + f(x)\ddt) + \sum_{i=1}^n \FDn(\k)\sdt
\cdot (\dx{i} + \dxsur{f}{i} \ddt)$.\\
{\bf Affirmation.} Cet id\'eal est l'annulateur de $f^s$ dans $\FDn(\k)\sdt$.
\item[1]
$I_1(f)= I(f) \cap \FDn(\k)[s]$.\\
Ainsi, $I_1(f)$ s'obtient \`a partir de $I(f)$ en \'eliminant la variable
(globale) $\ddt$.\\
{\bf Affirmation.} L'id\'eal $I_1(f)$ est l'annulateur de $f^s$ dans
$\FDn(\k)[s]$.
\item[2]
$\dps I_2(f) = I_1(f) + \FDn(\k)[s] \cdot f$.
\item[3]
$\dps J(f) = I_2(f) \cap \k[[x]][s]$.\\
Ainsi, $J(f)$ s'obtient en \'eliminant les variables (globales) $\dx{i}$.
\item[4]
\begin{eqnarray*}
\B(f) & = & J(f) \cap \k[s]\\
      & = & I_2(f) \cap \k[s].
\end{eqnarray*}
Ce dernier id\'eal s'obtient en ``\'eliminant'' les variables (locales)
$x_i$ dans l'id\'eal $J(f)$.\\
{\bf Affirmation.} L'id\'eal $\B(f)$ est l'id\'eal de Bernstein de $f$.
Son g\'en\'erateur unitaire qu'on note $b_f$ est le polyn\^ome de Bernstein
de $f$.
\end{description}

\begin{proof}[D\'emonstration des affirmations]
D\'emontrons la premi\`ere. Il est facile de voir que $I(f)$ est inclus
dans l'annulateur de $f^s$. Montrons l'inclusion inverse. Soit $P \in
\FDn(\k)\sdt$ s'annulant sur $f^s$. Modulo $I(f)$, on peut supposer que
$P$ appartient \`a $\k[[x]][\ddt]$. \'Ecrivons $P=\sum_\nu u_\nu(x)\ddt^\nu$.
On a alors $P\cdot f^s=\sum_\nu u_\nu(x) (-1)^\nu s (s-1) \cdots (s-\nu+1)
f^{-\nu} f^s=0$. Cette \'egalit\'e ayant lieu dans $\mathcal{L}$, on en
d\'eduit que les $u_\nu(x)$ sont nuls ce qui d\'emontre la premi\`ere
affirmation. La seconde \'etant triviale, voyons la troisi\`eme. Soit $c(s)
\in \k[s]$. C'est un polyn\^ome de Bernstein de $f$ si et seulement si il
existe $P\in \FDn(\k)[s]$ tel que $c(s) f^s=P \cdot f^{s+1}$, ou encore
$c(s)-Pf$ annule $f^s$, ce qui d'apr\`es la seconde affirmation est
\'equivalent \`a $c(s) -P f \in I_1(f)$ ou encore $c(s) \in I_1(f)
+\FDn(\k)[s] f$.
\end{proof}

D'apr\`es les r\'esultats de la section pr\'ec\'edente, nous savons calculer
des g\'en\'erateurs des id\'eaux $I_1(f)$, $I_2(f)$ et $J(f)$, ceci
en faisant un calcul de bases standard pour un ordre bien choisi. Le
probl\`eme est maintenant le suivant~: \'etant donn\'e $J\subset \k[[x]][s]$,
comment calculer le g\'en\'erateur unitaire $b$ de $J\cap \k[s]$~?
Nous supposerons ce $b$ non nul (ce qui est le cas dans notre situation).

Dans (\cite{oaku}, Algorithme 4.5), T.~Oaku a trait\'e la question suivante~:
soit $J\subset \k[x][s]$, comment calculer $\big(\k[[x]][s]\cdot J\big)
\cap \k[s]$~?

Nous allons utiliser le m\^eme algorithme en apportant une l\'eg\`ere
modification \`a la d\'emonstration. Pour les besoins du probl\`eme, nous
aurons besoin de travailler avec la cl\^oture alg\'ebrique $\bar{\k}$ de $\k$.
Cependant, nous verrons que si l'on sait \`a l'avance que $b(s)$ est \`a
racines dans $\k$ alors $\bar{\k}$ est inutile.

\begin{num}\label{algo:oakuformel}
\'Elimination des variables $x_i$: version formelle d'un algorithme de Oaku.
\end{num}

\begin{description}
\item[$(\alpha)$] Soit $b_0(s)$ le g\'en\'erateur unitaire de
$J(0,s)=\{g(0,s)|g(x,s) \in J\}$ (qui forme un id\'eal de
$\k[s]$).\\ Remarquons que $J(0,s)$ est engendr\'e par $\{g(0,s)| g\in
G\}$ si $G$ engendre $J$. Ainsi $b_0(s)$ s'obtient via un calcul de
pgcd ou bien un calcul de bases de Gr\"obner (ce que nous
utiliserons).\\ Remarque~: $b(s)$ est un multiple de $b_0(s)$ qui est
donc non nul.

\item[$(\beta)$] Soit $b_0(s)=(s-s_1)^{\mu_1} \cdots (s-s_m)^{\mu_m}$
la factorisation de $b_0(s)$ dans $\bar{\k}[s]$.\\ Par la remarque
pr\'ec\'edente, $b(s)$ s'\'ecrit $b(s)=p(s) (s-s_1)^{\nu_1} \cdots
(s-s_m)^{\nu_m}$ avec $p(s_i) \ne 0$ et $\nu_i \ge \mu_i$ pour
$i=1,\ldots ,m$.

\item[$(\gamma)$]
Soit $\bar{J} = \bar{\k}[[x]][s] \cdot J$.

\item[$(\delta)$] Pour $i=1,\ldots,m$, soit $l_i \in \N$ le plus petit
entier $l$ tel qu'il existe $h(x,s) \in \bar{\k}[[x]][s]$ avec $h(x,s)
(s-s_i)^l \in \bar{J}$ et $h(0,s_i) \ne 0$.\\ En consid\'erant $b(s)$,
on constate que de tels $l$ et $h(x,s)$ existent et que $l_i \le
\nu_i$. De plus, en faisant $(x,s)=(0,s_i)$, on voit que $l_i \ge \mu_i$.\\
Enfin remarquons que pour un $l$ donn\'e, un tel $h(x,s)$ se trouve dans
le quotient $\bar{J} : (s-s_i)^l$. Ainsi, on voit ais\'ement que $l\ge l_i$
si et seulement si n'importe quel syst\`eme de g\'en\'erateurs de $\bar{J}~:
(s-s_i)^l$ contient un \'el\'ement qui ne s'annule pas en
$(x,s)=(0,s_i)$.

\item[$(\varepsilon)$]
On pose $c(s) = (s-s_1)^{l_1} \cdots (s-s_m)^{l_m}$.
\end{description}

\begin{prop}
$b(s)$ est \'egal \`a $c(s)$.
\end{prop}

\begin{proof}
Posons $E=\bar{J} : c(s)$; cet id\'eal contient
$\bar{J}$. Consid\'erons le lieu des z\'eros de $E(0,s)\subset
\bar{\k}[s]$ dans $\bar{\k}$. On a alors~:
\[V(E(0,s)) \subset V(J(0,s))= V(b_0(s))=\{s_1,\ldots,s_m\}.\]
D'autre part, pour chaque $i=1,\ldots,m$, il existe $h_i(x,s) \in
\bar{\k}[[x]][s]$ tel que $h_i(x,s) (s-s_i)^{l_i} \in \bar{J}$ et
$h_i(0,s_i) \ne 0$. Ainsi, $\dps h_i(x,s)$ appartient \`a $E$ et ne
s'annule pas en $(x,s)=(0,s_i)$. En cons\'equence, $V(E(0,s))= \emptyset$.

Par le th\'eor\`eme des z\'eros de Hilbert, $1 \in E(0,s)$. Cela signifie
qu'il existe $e=e(x,s) \in E$ tel que $e(0,s)=1$. Quitte \`a multiplier
$e(x,s)$ par une unit\'e de $\bar{\k}[[x]]$, on peut supposer que $e
\in 1+\sum_{i=1}^n \bar{k}[[x]][s] \cdot (x_i s)$.

Maintenant, notons $d$ le degr\'e de $b(s) \in J \subset \bar{J}
\subset E$. Consid\'erons le $\bar{\k}[[x]]$-module $M=\bar{\k}[[x]]
\oplus \cdots \oplus \bar{\k}[[x]] s^d$ et posons $N=M \cap E$.

Montrons que pour tout entier $q$, $1$ appartient \`a $m^q M +N$, $m$ \'etant
l'id\'eal maximal de $\bar{\k}[[x]]$.

D'apr\`es ce qui pr\'ec\`ede, il existe $v\in m \bar{\k}[[x]][s] s$ et $e\in
E$ tel que $1=v+e$. En \'elevant \`a la puissance $q$, on obtient~:
\[1 \in m^q \bar{\k}[[x]][s]s + E.\]
Ecrivons~: $1=v_1 s+ v_2 s^2 + \cdots + v_N s^N + e'$ avec $v_i \in m^q$
et $e'\in E$. Pour chaque $k=N, N-1, \ldots, d+1$ (si $N\ge d+1$), on
divise $v_k s^k$ par $b(s)$ et on obtient $v_k s^k \in m^q s + \cdots
+ m^q s^{k-1} +E$. \`A la fin de ces divisions, on a
$1\in m^q s + \cdots + m^q s^d +E$.
Ainsi $1 \in m^q M + N$. Par le th\'eor\`eme d'intersection de Krull, $1 \in
N \subset E$, i.e. $c(s) \in \bar{J} \cap \bar{\k}[s]$.

Remarque~: le polyn\^ome $b(s)$ joue ici le r\^ole du polyn\^ome
$g(x,s)$ de \cite[Algo. 4.5]{oaku}. Afin d'arriver \`a $1\in
E$, T.~Oaku \emph{loc. cit.} a utilis\'e le th\'eor\`eme dit
d'extension \cite[Chap. 3, \S 6]{clo}.

Montrons que $c(s)$ engendre $\bar{J} \cap \bar{\k}[s]$. Soit $t(s)
\in \bar{J} \cap \bar{\k}[s]$. Ecrivons $t(s)=q(s) (s-s_1)^{u_1}
\cdots (s-s_m)^{u_m}$ avec $q(s_i) \ne 0$. Par d\'efinition de $l_i$,
on a $u_i \ge l_i$ ainsi $t(s)$ est multiple de $c(s)$.

Pour finir, montrons que $c(s)$ est dans $J \cap \k[s]$. Notons
$\pi: \bar{\k} \to \k$ la projection $\k \oplus S \to \k$ o\`u $S$ est
un suppl\'ementaire de $\k$. On l'\'etend \`a $\bar{\k}[[x]][s]$.
Comme $c(s)$ appartient \`a $\bar{J}$ et
que ce dernier est engendr\'e par $J$, on peut \'ecrire~:
$c(s)=\sum_i q_i(x,s) g_i(x,s)$
o\`u les $q_i$ sont dans $\bar{\k}[[x]][s]$ et les $g_i$ dans $J$
(donc dans $\k[[x]][s]$). Appliquons $\pi$ et remarquons que puisque
$g_i \in \k[[x]][s]$, on a $\pi(q_ig_i)=\pi(q_i) g_i$. Nous
obtenons que $\pi(c(s))$ appartient \`a $J \cap \k[s]$ donc \`a
$\bar{J} \cap \bar{\k}[s]$ et par cons\'equent est multiple de
$c(s)$. Or, puisque $c(s)$ est unitaire, il a m\^eme degr\'e que
$\pi(c(s))$. Ainsi $c(s)$ \'egale $\pi(c(s))$ et appartient bien \`a
$J \cap \k[s]$.

Nous savions que $c(s)$ divise $b(s)$. Maintenant, nous savons que $c(s)$
est multiple de $b(s)$ ce qui ach\`eve la preuve de la proposition.
\end{proof}

Dans la suite, lorsque nous utiliserons l'algorithme pr\'ec\'edent,
nous serons dans une situation o\`u l'on sait \`a l'avance que le
polyn\^ome de Bernstein est \`a racines rationnelles. Cela nous permet
de simplifier l'algorithme de la fa\c{c}on suivante.

\begin{rem}\label{rem:sanscloture}
Supposons que dans l'algorithme pr\'ec\'edent, les racines de $b(s)$ soient
dans $\k$ alors $\bar{\k}$ est inutile, plus pr\'ecis\'ement~:
\begin{itemize}
\item
Dans l'\'etape $(\beta)$, la factorisation se fait dans $\k[s]$.
\item
L'\'etape $(\gamma)$ peut \^etre saut\'ee.
\item
Enfin, dans l'\'etape $(\delta)$, il suffit de consid\'erer $J:(s-s_i)^l$,
i.e. $l_i$ est le plus petit $l$ tel que $J:(s-s_i)^l$ contienne un $h(x,s)
\in \k[[x]][s]$ ne s'annulant pas en $(x,s)=(0,s_i)$.
\end{itemize}
\end{rem}

\section{D\'emonstration de {\bf (i)} et {\bf (iii)} du th\'eor\`eme \ref{t}}

\`A partir d'ici, $\CC=\O_Y$. Consid\'erons les id\'eaux suivants:
\begin{description}
\item[0]
$\dps I=\FDn(\CC)\sdt\cdot (s + f(x)\ddt) + \sum_{i=1}^n \FDn(\CC)\sdt
\cdot (\dx{i} + \dxsur{f}{i} \ddt)$.
\item[1]
$I_1= \big( I+ \FDn(\QQ)\sdt \big) \cap \FDn(\CC)[s]$.
\item[2]\label{I_2}
$\dps I_2 = I_1 + \FDn(\CC)[s] \cdot f$.
\item[3]
$\dps J = \big( I_2+ \FDn(\QQ)[s] \big) \cap \CC[[x]][s]$.
\end{description}

Soit $\GG_0$ une base standard g\'en\'erique de $I$ sur $V(\QQ)$ pour un
ordre $\prec_0$ qui \'elimine la variable $\ddt$. Soit $h_0 \in \CC
\smallsetminus \QQ$ le produit des $\cp_{\prec_0}^\modQ(P)$ pour $P\in \GG_0$
et soit $\GG'_0$ le sous-ensemble de $\GG_0$ constitu\'e d'\'el\'ements
dont le terme privil\'egi\'e modulo $\QQ$ est ind\'ependant de $\ddt$.
De m\^eme, $\GG_2$ est une base standard g\'en\'erique de $I_2$ sur $V(\QQ)$
pour un ordre $\prec_2$ qui \'elimine les variables $\dx{i}$. On note $h_2
\in \CC \smallsetminus \QQ$ le produit des coefficients privil\'egi\'es
modulo $\QQ$ et on d\'efinit $\GG_2' \subset \GG_2$ comme le sous-ensemble
dont les \'el\'ements ont leur $\tp_{\prec_2}^\modQ$ ind\'ependant des
$\dx{i}$.

\begin{lem}
\begin{itemize}
\item[(0)]
Pour tout $\PP\in \spec(\CC)$, $(I)_\PP=I((f)_\PP)$.
\item[(1)]
Pour tout $\PP \in V(\QQ) \smallsetminus V(h_0)$, $(I_1)_\PP= I_1((f)_\PP)$.
\item[(2)]
Pour tout $\PP \in V(\QQ) \smallsetminus V(h_0)$, $(I_2)_\PP= I_2((f)_\PP)$.
\item[(3)]
Pour tout $\PP \in V(\QQ) \smallsetminus V(h_0 h_2)$, $(J)_\PP=J((f)_\PP)$.
\end{itemize}
\end{lem}
Rappelons que les notations $I((f)_\PP)$, $I_1((f)_\PP)$, etc, sont
celles introduites dans la construction algorithmique formelle
donn\'ee page \pageref{debutalgo} (ici on applique la construction \`a
$(f)_\PP \in \FF(\PP)[[x]]$).
\begin{proof}
L'assertion (0) est triviale. L'assertion (1) d\'ecoule directement du
corollaire \ref{cor:eliminegen} et de (0). L'assertion (2) est
une cons\'equence directe de (1). La (3) d\'ecoule du corollaire
\ref{cor:eliminegen} et de (2).
\end{proof}

Nous en sommes \`a la fin de l'\'etape $3$. Pour poursuivre, nous avons
besoin de quelques r\'esultats suppl\'ementaires.

\subsection{R\'esultats pr\'eparatoires et d\'ebut de la fin}

Dans la suite, nous aurons besoins de calculer ``g\'en\'eriquement''
les quotients du type $J:u$ o\`u $J\subset \k[[x]][s]$ et $u \in
\k[s]$. Rappelons comment les calculer dans le cas absolu (i.e. non
param\'etrique).

\begin{lem}[\cite{clo}, chap. 4, \S{}3 et \S{}4]\label{lem:clo}
Soit $\zeta$ une nouvelle variable, alors
\[J \cap (\k[[x]][s]\cdot u) = \Big( \k[[x]][s][\zeta] \cdot \zeta \cdot J+
\k[[x]][s][\zeta]\cdot (1-\zeta)\cdot u \Big) \cap \k[[x]][s]\]
et si $G=\{g_1, \ldots, g_r\}$ est un syst\`eme de g\'en\'erateurs de
$J \cap \k[[x]][s]\cdot u$ alors $G/u=\{g_1/u, \ldots, g_r/u\}$
engendre $J:u$.
\end{lem}
Ainsi, le calcul se r\'esume en la simple \'elimination d'une variable
globale.

Pour finir cette sous-section, voici un r\'esultat indispensable pour
continuer. C'est lui qui nous permettra, via la remarque
\ref{rem:sanscloture}, de d\'emontrer le point {\bf (i)} du th\'eor\`eme
\ref{t}.
\begin{lem}[de rationalit\'e]\label{lem:ratio}
Soit $p \in \O_Y[s]$ dont on note $\cp(p)\in \O_Y$ le coefficient du mon\^ome
de plus haut de degr\'e. Supposons qu'il existe un ouvert $W$ de Zariski
de $V_m(\QQ) \subset Y$ tel que pour tout $y\in W$, $(p/\cp(p))_{|y_0}$
(soit bien d\'efini) et appartienne \`a $\Q[s]$ alors il existe $q\in \Q[s]$
unitaire tel que $p-\cp(p) q \in \QQ[s]$.
\end{lem}

\begin{proof}
\'Ecrivons $p=\sum_{i=1}^N c_i(y) s^i$ et $c_N=\cp(p)$. Pour chaque $i$,
consid\'erons l'application $D_i: y\in W \mapsto \frac{c_i(y)}{c_N(y)} \in
\C$. Son image $D_i(W)$ est un constructible de $\C$
(voir \cite[Th. 3.16]{harris}),
or par hypoth\`ese $D_i(W) \subset \Q$ donc $D_i(W)$
est une r\'eunion finie de points rationnels. Par l'irr\'eductibilit\'e de
$V_m(\QQ)$, $D_i(W)$ est un singleton. Par cons\'equent, il existe
$q_i \in \Q$ tel que $c_i(y) -q_i c_N(y) \in \QQ$. Le polyn\^ome
$q =\sum_i q_i s^i$ est le polyn\^ome que l'on cherchait.
%
%
\end{proof}

\subsection{Les \'etapes $(\alpha)$ et $(\beta)$}

Reprenons la d\'emonstration du th\'eor\`eme \ref{t}. Nous en
\'etions \`a la fin de l'\'etape {\bf 3} o\`u nous avions construit
l'id\'eal $J \subset \CC[[x]][s]$.

\begin{lem}
Pour tout $\PP \in V(\QQ) \smallsetminus V(h_0h_2)$,
\[\big( J(0,s)  \big)_\PP = (J)_\PP(0,s) = J( (f)_\PP)(0,s).\]
\end{lem}

La premi\`ere \'egalit\'e est triviale en utilisant les d\'efinitions
et la deux\`eme d\'ecoule directement de l'\'etape {\bf 3}.

Soit maintenant $\GG_3$ une base standard g\'en\'erique de
$J(0,s)$ sur $V(\QQ)$ relativement \`a l'ordre usuel de $\N$. Soit
$h_3$ le produit des coefficients privil\'egi\'e modulo $\QQ$ et soit
$\tilde{b_0}$ l'\'el\'ement de $\GG_3$ dont l'exposant privil\'egi\'e
modulo $\QQ$ est le plus petit (ou dit plus simplement, dont le degr\'e
en $s$ modulo $\QQ$ est le plus petit).

\begin{lem}
Pour tout $\PP \in V(\QQ) \smallsetminus V(h_0h_2h_3)$,
$(\tilde{b_0})_\PP$ engendre $\big( J(0,s) \big)_\PP$.
\end{lem}

Ce lemme d\'ecoule du corollaire \ref{cor:eliminegen}.
Afin de poursuivre dans de bonnes conditions, nous avons besoin du

\begin{lem}
Il existe $b_0(s)$ unitaire et \`a racines rationnelles tel que
\[\tilde{b_0}(s) - \cp(\tilde{b_0}(s)) \cdot b_0(s) \in \QQ[s].\]
\end{lem}

\begin{proof}
Par l'algorithme \ref{algo:oakuformel}, nous savons que pour tout
$y_0\in Y$, le polyn\^ome de Bernstein de $f(x,y_0)$ (qui est \`a racines
dans $\Q$) a les m\^emes racines que le g\'en\'erateur de
$J(f(x,y_0))(0,s)$. Ainsi d'apr\`es ce qui pr\'ec\`ede, pour tout
$y_0$ dans un ouvert de Zariski de $V_m(\QQ)$, $(\tilde{b_0})_{|y_0}$ est
\`a racines dans $\Q$. On peut donc appliquer le lemme de rationalit\'e
\ref{lem:ratio} ce qui nous fournit $b_0 \in \Q[s]$ unitaire v\'erifiant
la relation $\tilde{b_0}(s) - \cp(\tilde{b_0}(s)) \cdot b_0(s) \in \QQ[s]$.
En sp\'ecialisant encore dans un ouvert de $V_m(\QQ)$ on montre
que ce $b_0$ est \`a son tour \`a racines rationnelles.
\end{proof}

Notons qu'\`a ce stade de la d\'emonstration, nous savons que pour un
$\PP$ g\'en\'erique dans $V(\QQ)$, les racines du polyn\^ome de
Bernstein de $(f)_\PP$ (qui sont celles de $b_0(s)$) sont rationnelles
et constantes. En particulier, {\bf le point {\bf (i)} du th\'eor\`eme
\ref{t} est acquis}.

\subsection{Fin du parcours~: l'\'etape $(\delta)$}

Pour le moment, nous savons que pour tout $\PP \in V(\QQ) \smallsetminus
V(h_0h_2h_3)$, les polyn\^omes suivant sont \'egaux~:
\begin{itemize}
\item
le polyn\^ome $b_0$ (\`a racines rationnelles) obtenu dans le lemme
pr\'ec\'edent,
\item
le polyn\^ome qu'on note $b_0((f)_\PP)$ et qui est celui qu'on obtient
\`a l'\'etape $(\alpha)$ de l'algorithme \ref{algo:oakuformel} appliqu\'e
\`a $(f)_\PP$.
\end{itemize}
Au passage, introduisons quelques autres notations. Soient
$s_1,\ldots,s_m \in \Q$ les racines de $b_0$. Pour chaque $i$, et
chaque $\PP$ comme au dessus, on note $\mu_i(\PP)$ et $l_i(\PP)$ les
entiers obtenus dans l'algorithme \ref{algo:oakuformel} appliqu\'e \`a
$(f)_\PP$.

On sait d\'ej\`a que les $\mu_i(\PP)$ sont g\'en\'eriquement constants et
tout ce qui nous reste \`a faire, c'est de montrer qu'en excluant une
nouvelle hypersurface de $V(\QQ)$, les $l_i(\PP)$ sont aussi constants.

Fixons $i\in \{1,\ldots,m\}$ et posons $\mu_i=\mu_i(\QQ)$ et
$l_i=l_i(\QQ)$. Maintenant, pour tout entier $l$ avec $\mu_i \le l
\le l_i$, consid\'erons l'id\'eal $J(i,l)$ dans $ \CC[[x]][s][\zeta]$~:
\[ J(i,l) = \CC[[x]][s][\zeta] \cdot \zeta \cdot J + \CC[[x]][s][\zeta]
\cdot (1-\zeta)\cdot (s-s_i)^l. \]

Pour tout $\mu_i \le l \le l_i$, soit $\GG(i,l)$ une base standard
g\'en\'erique de $J(i,l)$ sur $V(\QQ)$ relativement \`a un ordre qui
\'elimine la variable $\zeta$. Notons $h_i'$ le produit des coefficients
privil\'egi\'es modulo $\QQ$ des \'el\'ements des $\GG(i,l)$ pour $l=\mu_i,
\ldots, l_i$. Finalement soit $\GG'(i,l)$ les \'el\'ements dont le terme
privil\'egi\'e modulo $\QQ$ est ind\'ependant de $\zeta$.

En utilisant le lemme \ref{lem:clo}, on a: pour tout $\PP \in V(\QQ)
\smallsetminus V(h_0h_2h_3h_i')$, $(\GG'(i,l))_\PP$ engendre
\begin{eqnarray*}
(J(i,l))_\PP \cap \FF(\PP)[[x]][s] & = & \Big( \zeta \cdot (J)_\PP +
  (1-\zeta)(s-s_i)^l \Big) \cap \FF(\PP)[[x]][s] \\
& = & (J)_\PP \cap \big( \FF(\PP)[[x]][s] \cdot (s-s_i)^l \big).
\end{eqnarray*}

Ainsi, pour les m\^emes $\PP$, $\dps \frac{(\GG'(i,l))_\PP}{(s-s_i)^l}$
engendre $(J)_\PP : (s-s_i)^l$.

\begin{rem}
Soit $P \in \GG'(i,l)$ alors le diagramme de Newton de
$(P)_\QQ$ est \'egal \`a $\DN^\modQ(P)$.
Or $\frac{1}{(s-s_i)^l}(P)_\QQ$ est dans $\FF(\QQ)[[x]][s]$
(i.e. n'a pas de p\^ole en $s=s_i$) donc il existe un unique couple
$(P^1,P^2)$ avec $P^1 \in (\CC \smallsetminus \QQ)[[x]][s] \cdot (s-s_i)^l$
et $P^2 \in \QQ[[x]][s][\zeta]$ tel que $P=P^1+P^2$.
Ici, $(\CC \smallsetminus \QQ)[[x]][s]$ est le sous-ensemble de $\CC[[x]][s]$
dont les \'el\'ements ont leur coefficients dans $\CC \smallsetminus \QQ$.
\end{rem}
\noindent
$\bullet$ Montrons que, pour $\PP \in V(\QQ) \smallsetminus
V(h_0h_2h_3h_i')$, $l_i(\PP) \ge l_i=l_i(\QQ)$.\\
Par l'absurde, soit $l< l_i$ pour lequel il existe un \'el\'ement dans
$(J)_\PP : (s-s_i)^l$ qui ne s'annule pas en $(x,s)=(0,s_i)$.
Cela signifie qu'il existe $P$ dans $\GG'(i,l)$ tel que
$\dps \Big( \frac{P}{(s-s_i)^l}\Big)_\PP (0,s_i) \ne 0$ (voir les remarques
\`a l'\'etape $(\delta)$ de l'algorithme \ref{algo:oakuformel}).
Or $\QQ \subset \PP$ donc $\dps \Big( \frac{P}{(s-s_i)^l} \Big)_\QQ
(0,s_i) \ne 0$. Cela contredit la propri\'et\'e de minimalit\'e de
$l_i(\QQ)$.

\noindent
$\bullet$ Il nous reste \`a montrer l'in\'egalit\'e
r\'eciproque. Pour cela, nous allons exclure une derni\`ere
hypersurface.

Par d\'efinition de $l_i=l_i(\QQ)$, il existe $P$ dans $\GG'(i,l_i)$ tel
que\\$\dps \Big( \frac{P}{(s-s_i)^l} \Big)_\QQ (0,s_i) \ne 0$.
Pour un tel $P$, notons
\[h_i''= \frac{P^1}{(s-s_i)^{l_i}}(0,s_i) \in \CC \smallsetminus \QQ.\]
Voir la remarque ci-dessus pour la d\'efinition de $P^1$.
Pour tout $\PP \in V(\QQ) \smallsetminus V(h_0h_2h_3h_i' h_i'')$,
$\dps \Big( \frac{P}{(s-s_i)^l} \Big)_\PP (0,s_i) \ne 0$, ce qui
implique $l_i(\PP) \le l_i$.\\ 
\underline{Bilan~:}\\
Si on note $h'$ le produit de $h_0h_2h_3$, des $h_i'$ et des $h_i''$,
alors pour tout $\PP \in V(\QQ) \smallsetminus V(h')$,
le polyn\^ome de Bernstein (formel) de $(f)_\PP$ est constant.
Le point {\bf (iii)} du th\'eor\`eme \ref{t} en d\'ecoule puisque
le polyn\^ome de Bernstein analytique de $f_y$ est \'egal \`a son
homologue formel \cite{bmai90}.

\section{L'assertion {\bf (ii)} du th\'eor\`eme \ref{t}}

Si $c(s) \in \C[s]$ satisfait (\ref{eq:faible}) alors une sp\'ecialisation
en un $y$ g\'en\'erique de $V_m(\QQ)$ nous dit que $c(s)$ est un multiple
de $b(s)$ (ceci gr\^ace \`a {\bf (iii)}). Ainsi pour d\'emontrer {\bf (ii)},
il suffit de montrer que $b(s)$ satisfait (\ref{eq:faible}).

Soit $\GG$ une base standard g\'en\'erique de $I_2$ (cf. page \pageref{I_2})
sur $V(\QQ)$ pour
un ordre $\prec$ quelconque (on demande juste que le th\'eor\`eme de
division dans $\FDn(\k)[s]$ par rapport \`a cet ordre existe).
Notons $c$ le produit des coefficients privil\'egi\'es modulo $\QQ$ de $\GG$.
Par construction, on sait que $b=(b)_\QQ$ appartient \`a $(I_2)_\QQ$.
On peut donc effectuer la division modulo $\QQ$ de $b$ par $\GG$ dans
$\FDn(\O_Y)[s]$, division dont le reste est nul modulo $\QQ$ par le
corollaire \ref{cor:restemodQ}.
Ainsi~: $b(s) \in I_2 + \FDn(\QQ)[s]$. En appliquant $b(s)$ \`a $f^s$, on
obtient une \'equation fonctionnelle
\[b(s) f^s = P_0(s) f^{s+1} + P_1(s,\ddt) f^s\]
o\`u $P_0(s)$ appartient \`a $\FDn(\O_Y[c^{-1}])[s]$ et $P_1(s,\ddt)$
\`a $\FDn(\QQ[c^{-1}])\sdt$.

\`A partir de cette \'equation, nous allons faire un passage du formel \`a
l'analytique en nous inspirant de \cite{bmai90}.

On consid\`ere $\FDn(\O_Y[c^{-1}])[s]f^{s+1}$ et
$\FDn(\QQ[c^{-1}])\sdt f^s$ dans la somme directe suivante
\[\bigoplus_{l \ge 0} \O_Y[c^{-1}][[x]] \xi_l\subset
\O_Y[c^{-1}][[x]][1/f,s] f^s\]
o\`u $\xi_0=f^s$ et $\xi_l=(s-l+1)\cdots s f^{s-l}$ pour $l\ge 1$.

L'action de $\FDn(\O_Y[c^{-1}])\sdt$ se r\'esume \`a~: pour $u \in
\O_Y[c^{-1}][[x]]$,
\begin{itemize}
\item
$\dxi \cdot u \xi_l = \dxsur{u}{i} \xi_l + u \dxsur{f}{i} \xi_{l+1}$,
\item
$\ddt \cdot u \xi_l= -u \xi_{l+1}$,
\item
$s \cdot u \xi_l= lu \xi_l + u f \xi_{l+1}$.
\end{itemize}

Maintenant, soit $d$ le maximum des entiers $\deg_s(b(s))$,
$\deg_{\ddx,s}(P_0(s))$ et $\deg_{\ddx,s,\ddt} (P_1(s,\ddt))$.  Au vu des
trois identit\'es ci-dessus, l'\'equation du d\'ebut a lieu dans $\dps
\bigoplus_{0\le l\le d} \O_Y[c^{-1}][[x]] \xi_l$.
Remarquons d\`es \`a pr\'esent que~:
\[b(s) f^s \in \bigoplus_{0\le l\le d} \O_Z \xi_l.\]
Ecrivons $\dps P_0(s)=\sum V_0(\beta,k) \ddx^\beta s^k$ avec
$V_0(\beta,k) \in \O_Y[c^{-1}][[x]]$. Pour chaque $(\beta,k)$,
\'ecrivons $\dps \ddx^\beta s^k f^{s+1} = \sum_{0\le l\le d}
\phi_0(\beta,k,l)  \xi_l$.
Remarquons que les $\phi_0(\beta,k,l)$ sont dans $\O_Z$.

De m\^eme, on \'ecrit $\dps P_1(s,\ddt) =\sum V_1(\beta,k,\nu)
\ddx^\beta s^k \ddt^\nu$ avec $V_1(\beta,k,\nu)$ dans
$\QQ[c^{-1}][[x]]$ et pour chaque $(\beta,k,\nu)$, $\dps \ddx^\beta
s^k \ddt^\nu f^s=\sum_{0\le l\le d} \phi_1(\beta,k,\nu,l) \xi_l$.
Les $\phi_1(\beta,k,\nu,l)$ sont aussi dans $\O_Z$.\\

Maintenant, si on note $K_0$ le nombre de $V_{\beta,k}$, $V_0 \in
\big( \O_Y[c^{-1}][[x]] \big)^{K_0}$ le vecteur form\'e de ces
derniers et $\Phi_0$ la matrice $K_0\times d$ form\'ee des
$\phi_0(\beta,k,l)$ alors on peut repr\'esenter $P_0(s) f^{s+1}$ sous la
forme $\Phi_0 (V_0)$.

On d\'efinit de la m\^eme mani\`ere $K_1$, $V_1 \in \big(
\QQ[c^{-1}][[x]] \big)^{K_1}$ et $\Phi_1$ (matrice de taille
$K_1\times d$) et on a~: $P_1(s,\ddt) f^s$ s'identifie \`a $\Phi_1(V_1)$.

Notons que les matrices $\Phi_0$ et $\Phi_1$ sont \`a coefficients
dans $\O_Z \subset \O_Z[c^{-1}]$ ce qui nous permet de
d\'efinir les applications $\O_Z[c^{-1}]$-lin\'eaires~:
\begin{itemize}
\item
$\Phi'_0 : \big( \O_Z[c^{-1}] \big)^{K_0} \to \big(
\O_Z[c^{-1}] \big)^d$
\item
$\Phi'_1 : \big( \O_Z[c^{-1}]\cdot \QQ \big)^{K_0} \to \big(
\O_Z[c^{-1}] \big)^d$
\end{itemize}
donn\'ees respectivement par les matrices en question.

Enfin, notons $W \in \big( \O_Z \big)^d$ le vecteur
repr\'esentant $b(s)f^s$ dans la base $\xi_0, \ldots,
\xi_l$. Remarquons qu'on peut voir $W$ dans $\big( \O_Z[c^{-1}]
\big)^d$.

Maintenant, soit $q\ge 0$ un entier quelconque.
Consid\'erons la troncature de $V_0$ \`a l'ordre $q$, i.e. \'ecrivons
$V_0=V'_0 + V''_0$ o\`u $V'_0 \in \big( \O_Y[c^{-1}][x] \big)^{K_0}$ est
de degr\'e (en $x$) plus petit que $q$, et $V''_0$ est
dans $m^q \cdot \big( \O_Y[c^{-1}][[x]] \big)^{K_0}$ o\`u $m$
d\'esigne l'id\'eal de $\O_Y[c^{-1}][[x]]$ engendr\'e par les $x_i$.
Le point cl\'e est que $V'_0$ est dans $\O_Z[c^{-1}]$.

Faisons de m\^eme pour $V_1=V'_1+ V''_1$ avec $V'_1 \in \big(
\QQ[c^{-1}][x] \big)^{K_0}$ et $V''_1 \in m^q \cdot \big(
\O_Y[c^{-1}][[x]] \big)^{K_0}$.
On a que $V'_1$ appartient \`a $\O_Z[c^{-1}] \cdot \QQ$.\\

L'\'equation fonctionnelle de d\'epart se traduit alors par les
\'egalit\'es
\begin{eqnarray*}
W & = & \Phi_0 \cdot V_0 + \Phi_1 \cdot V_1\\ & = & \Phi'_0 \cdot V'_0
  + \Phi'_1 \cdot V'_1 + \Phi_0 \cdot V''_0 + \Phi_1 \cdot V''_1.
\end{eqnarray*}
Or $W$ et $\Phi'_0 \cdot V'_0 + \Phi'_1 \cdot V'_1$ sont dans $\big(
\O_Z[c^{-1}]\big)^d$ donc $\Phi_0 \cdot V''_0 + \Phi_1 \cdot
V''_1$ appartient \`a $\big( \O_Z[c^{-1}]\big)^d$ et sa valuation en $x$
est au moins $q$.

Si on note $m'$ l'id\'eal de $\O_Z[c^{-1}]$ engendr\'e par les
$x_i$, on obtient
\[W \in \mathrm{Im}(\Phi'_0) + \mathrm{Im}(\Phi'_1) +
(m')^q\cdot \big( \O_Z[c^{-1}]\big)^d.\]

Rappelons que $Z$ est un polydisque compact et donc que $\O_Z[c^{-1}]$
est noeth\'erien. Appliquons le th\'eor\`eme d'intersection de Krull
(voir \cite[page 150]{eisenbud}).
Il existe $p \in m'$ tel que
\[(1+p) W\in \mathrm{Im}(\Phi'_0) + \mathrm{Im}(\Phi'_1).\]
En remontant la construction, cette \'egalit\'e se traduit par:
\[(1+p)b(s) f^s \in \O_Z[c^{-1}]\langle \ddx \rangle [s] \cdot f^{s+1} +
\big(\O_Z[c^{-1}]\cdot \QQ \big) \langle \ddx \rangle \sdt \cdot f^s.\]
On multiplie alors par une puissance $k$ assez grande de $c$ pour chasser
les d\'enominateurs et l'on pose $h=c^k(1+p)$, ce qui fournit la relation
cherch\'ee.


\begin{thebibliography}{99}


\bibitem[ACG01]{acg01}
A. Assi, F. J. Castro-Jim\'enez, M. Granger,
\emph{The analytic standard fan of a $\mathcal{D}$-module},
J. Pure Appl. Algebra  164  (2001),  no. 1-2, 3--21.




\bibitem[Bah03a]{these}
R. Bahloul,
\emph{Contributions \`a l'\'etude des id\'eaux de Bernstein-Sato d'un point de vue constructif},
th\`ese de doctorat, Universit\'e d'Angers, 2003.
http://math.univ-anger.fr/$\sim$bahloul/




\bibitem[Bah03b]{PJA}
R. Bahloul,
\emph{Global generic Bernstein-Sato polynomial on an irreducible affine
scheme},
Proc. Japan Acad. Ser. A Math. Sci.  79  (2003),  no. 9, 146--149.


\bibitem[Bah04]{genSB}
R. Bahloul,
\emph{Generic and comprehensive standard bases},
preprint math.AC/0410220 (2004). Soumis.


\bibitem[Bah05]{Comp}
R. Bahloul,
\emph{D\'emonstration constructive de l'existence de polyn\^omes de
Bernstein-Sato pour plusieurs fonctions analytiques},
Compositio Math. 141  (2005),  no. 1, 175--191.


\bibitem[Ber72]{bernstein}
I. N. Bernstein,
\emph{Analytic continuation of generalized functions with respect
to a parameter},
Funkcional. Anal. i Prilo\v zen.  6  (1972), no. 4, 26--40.




\bibitem[Bio96a]{biosca}
H.~Biosca,
\emph{Sur l'existence de polyn\^omes de Bernstein g\'en\'eriques
associ\'es \`a une application analytique},
C. R. Acad. Sci. Paris S\'er. I Math.  322  (1996),  no. 7, 659--662.


\bibitem[Bio96b]{bioscaT}
H.~Biosca,
\emph{Polyn\^omes de Bernstein g\'en\'eriques et relatifs associ\'es
\`a une application analytique},
th\`ese de doctorat, Nice Sophia-Antipolis, 1996.


\bibitem[Bj\"o73]{bjorkP}
J. E. Bj\"ork,
\emph{Dimensions of modules over algebras of differential operators}, 
Fonctions analytiques de plusieurs variables et analyse complexe (Colloq.
Internat. CNRS, No. 208, Paris, 1972),  pp. 6--11.
``Agora Mathematica'', No. 1, Gauthier-Villars, Paris, 1974. 


\bibitem[Bj\"o79]{bjorkB}
J. E. Bj\"ork,
\emph{Rings of differential operators},
North-Holland Math. Library, 1979.


\bibitem[Bri]{briancon}
J. Brian\c{c}on,
\emph{Passage du local au global},
notes manuscrites.


\bibitem[BGM92]{bgm}
J.~Brian\c{c}on, F.~Geandier, Ph.~Maisonobe,
\emph{D\'eformation d'une singularit\'e isol\'ee d'hypersurface et
polyn\^omes de Bernstein},
Bull. Soc. Math. France  120  (1992),  no. 1, 15--49.


\bibitem[BGMM89]{bgmm}
J. Brian\c{c}on, M. Granger, Ph. Maisonobe, M. Miniconi,
\emph{Algorithme de calcul du polyn\^ome de Bernstein:
cas non d\'eg\'en\'er\'e},
Ann. Inst. Fourier (Grenoble)  39  (1989),  no. 3, 553--610.


\bibitem[BLM91]{blm}
J.~Brian\c{c}on, Y. Laurent, Ph.~Maisonobe,
\emph{Sur les modules diff\'erentiels holonomes r\'eguliers, coh\'erents
relativement \`a une projection},
C. R. Acad. Sci. Paris S\'er. I Math.  313  (1991),  no. 5, 285--288.


\bibitem[BM90]{bmai90}
J.~Brian\c{c}on, Ph.~Maisonobe,
\emph{Examen de passage du local au global pour les polyn\^omes de
Bernstein-Sato},
notes non publi\'ees, 1990.


\bibitem[BM02]{bmai02}
J.~Brian\c{c}on, Ph.~Maisonobe,
\emph{Remarques sur l'id\'eal de Bernstein associ\'e \`a des polyn\^omes},
preprint no. 650, Nice Sophia-Antipolis, 2002.






\bibitem[Buc70]{buchberger}
B. Buchberger,
\emph{Ein algorithmisches Kriterium f\"ur die L\"osbarkeit eines
algebraischen Gleichungssystems},
Aequationes Math.  4  (1970), 374--383.


\bibitem[Cas86]{cassou86}
P. Cassou-Nogu\`es,
\emph{Racines de polyn\^omes de Bernstein},
Ann. Inst. Fourier (Grenoble) 36 (1986), no. 4, 1--30.


\bibitem[Cas87]{cassou87}
P. Cassou-Nogu\`es,
\emph{\'Etude du comportement du polyn\^ome de Bernstein lors d'une
  d\'eformation \`a $\mu$-constant de $X\sp a+Y\sp b$ avec
  $(a,b)=1$},
Compositio Math. 63 (1987), no. 3, 291--313.


\bibitem[Cas88]{cassou88}
P. Cassou-Nogu\`es,
\emph{Polyn\^ome de Bernstein g\'en\'erique},
Abh. Math. Sem. Univ. Hamburg 58 (1988), 103--123.


\bibitem[CG04]{cg}
F.~J.~Castro-Jim\'enez, M.~Granger,
\emph{Explicit calculations in rings of differential operators},
\'El\'ements de la th\'eorie des syst\`emes diff\'erentiels g\'eom\'etriques,
89--128, S\'emin. Congr., 8, Soc. Math. France, Paris, 2004.


\bibitem[CLO92]{clo}
D.~Cox, J.~Little, D.~O'Shea,
\emph{Ideals, varieties, and algorithms},
Undergraduate Texts in Mathematics, Springer-Verlag, New York, 1992.
Seconde \'edition 1997.


\bibitem[Eis95]{eisenbud}
D. Eisenbud,
\emph{Commutative algebra with a view toward algebraic geometry},
Graduate Texts in Mathematics 150, Springer, New York, 1995.


\bibitem[Fri67]{frisch}
J. Frisch,
\emph{Points de platitude d'un espace analytique},
Invent. Math.  4  (1967), 118--138.


\bibitem[Gea89]{geandierCRAS}
F. Geandier,
\emph{Polyn\^omes de Bernstein et d\'eformations \`a nombre de Milnor
constant},
C. R. Acad. Sci. Paris S\'er. I Math.  309  (1989),  no. 13, 831--834.


\bibitem[Gea91]{geandierComp}
F.~Geandier,
\emph{D\'eformations \`a nombre de Milnor constant: quelques r\'esultats
sur les polyn\^omes de Bernstein},
Compositio Math.  77  (1991),  no. 2, 131--163.




\bibitem[Gyo93]{gyoja}
A.~Gyoja,
\emph{Bernstein-Sato's polynomial for several analytic functions},
J. Math. Kyoto Univ.  33  (1993),  no. 2, 399--411.


\bibitem[Har92]{harris}
J. Harris,
\emph{Algebraic Geometry. A first course},
Graduate Texts in Mathematics 133, Springer-Verlag, New York, 1992.




\bibitem[Kas76]{kashiwara}
M. Kashiwara,
\emph{$B$-functions and holonomic systems. Rationality of roots of
$B$-functions},
Invent. Math.  38  (1976/77), no. 1, 33--53.


\bibitem[L\^e73]{le}
D. T. L\^e,
\emph{Topologie des singularit\'es des hypersurfaces complexes},
Singularit\'es \`a Carg\`ese,  pp. 171--182.
Asterisque, Nos. 7 et 8, Soc. Math. France, Paris, 1973.


\bibitem[LR76]{le-r}
D. T. L\^e, C. P. Ramanujam,
\emph{The invariance of Milnor's number implies the invariance of
the topological type},
Amer. J. Math.  98  (1976), no. 1, 67--78.


\bibitem[Ley01]{leykin}
A. Leykin,
\emph{Constructibility of the set of polynomials with a fixed Bernstein-Sato
Polynomial: an algorithmic approach},
J. Symbolic Comput.  32  (2001),  no. 6, 663--675.


\bibitem[Mal74]{malgrange}
B. Malgrange,
\emph{Le polyn\^ome de Bernstein d'une singularit\'e isol\'ee},
Lecture Notes in Math. 459, 98--119, Springer, Berlin, 1975.





\bibitem[Oak97a]{oakuJPAA}
T. Oaku,
\emph{Algorithms for the $b$-function and $D$-modules associated with
a polynomial},
J. Pure Appl. Algebra  117/118  (1997), 495--518.


\bibitem[Oak97b]{oaku}
T. Oaku,
\emph{An algorithm of computing $b$-functions},
Duke Math. J. 87 (1997), no. 1, 115--132.




\bibitem[Sab87a]{sabbah1}
C. Sabbah,
\emph{Proximit\'e \'evanescente I. La structure polaire d'un
$\mathcal{D}$-Module},
Appendice en collaboration avec F.~J. Castro-Jim\'enez,
Compositio Math.  62  (1987),  no. 3, 283--328.


\bibitem[Sab87b]{sabbah2}
C. Sabbah,
\emph{Proximit\'e \'evanescente II. \'Equations fonctionnelles pour
plusieurs fonctions analytiques},
Compositio Math.  64  (1987),  no. 2, 213--241.


\bibitem[SS72]{sato}
M. Sato, T. Shintani,
\emph{On zeta functions associated with prehomogeneous vector spaces},
Proc. Nat. Acad. Sci. U.S.A.  69  (1972), 1081--1082;
Ann. of Math. (2)  100  (1974), 131--170.




\bibitem[Yan78]{yano}
T. Yano,
\emph{On the theory of $b$-functions},
Publ. Res. Inst. Math. Sci.  14  (1978), no. 1, 111--202.


\end{thebibliography}
\end{document}